\newcommand{\xbG}{\Gamma}

\newcommand{\xba}{\alpha}
\newcommand{\xbb}{\beta}

\newcommand{\xbe}{\in}
\newcommand{\xbf}{\phi}
\newcommand{\xbg}{\gamma}

\newcommand{\xbm}{\mu}

\newcommand{\xbq}{\psi}
\newcommand{\xbr}{\rho}
\newcommand{\xbs}{\sigma}
\newcommand{\xbt}{\tau}

\newcommand{\xCB}{A}

\newcommand{\xCN}{\neg}
\newcommand{\xCQ}{\emptyset}

\newcommand{\xcN}{\hspace{0.2em}\not\sim\hspace{-0.9em}\mid\hspace{0.8em}}

\newcommand{\xcP}{\not\rightarrow}

\newcommand{\xcS}{\bigcap}
\newcommand{\xcT}{\bot}

\newcommand{\xcc}{\subseteq}

\newcommand{\xch}{\Rightarrow}

\newcommand{\xcj}{\Leftrightarrow}
\newcommand{\xck}{\leq}
\newcommand{\xcl}{\vdash}
\newcommand{\xcm}{\models}
\newcommand{\xcn}{\hspace{0.2em}\sim\hspace{-0.9em}\mid\hspace{0.58em}}

\newcommand{\xco}{\vee}
\newcommand{\xcp}{\rightarrow}

\newcommand{\xcr}{\leftrightarrow}
\newcommand{\xcs}{\cap}
\newcommand{\xcu}{\wedge}
\newcommand{\xcv}{\cup}

\newcommand{\xcz}{\Box}

\newcommand{\xDH}{\item }

\newcommand{\xDM}{\circ}

\newcommand{\xdC}{\mbox{\boldmath$C$}}

\newcommand{\xdf}{{\cal F}}

\newcommand{\xdi}{{\cal I}}

\newcommand{\xdl}{{\cal L}}

\newcommand{\xdp}{{\cal P}}

\newcommand{\xdy}{{\cal Y}}

\newcommand{\xEH}{ & }
\newcommand{\xEI}{\begin{itemize}}
\newcommand{\xEJ}{\end{itemize}}
\newcommand{\xEP}{ \\ }

\newcommand{\xEc}{\not<}
\newcommand{\xEd}{\neq}
\newcommand{\xEh}{\begin{enumerate}}
\newcommand{\xEj}{\end{enumerate}}

\newcommand{\xFB}{\cdots}

\newcommand{\Xl}{\ldots}

\newcommand{\ol}{\overline}

\newcommand{\bl}{\begin{lemma} \rm}
\newcommand{\el}{\end{lemma}}
\newcommand{\br}{\begin{remark} \rm}
\newcommand{\er}{\end{remark}}
\newcommand{\be}{\begin{example} \rm}
\newcommand{\ee}{\end{example}}
\newcommand{\bco}{\begin{corollary} \rm}
\newcommand{\eco}{\end{corollary}}
\newcommand{\bc}{\begin{claim} \rm}
\newcommand{\ec}{\end{claim}}
\newcommand{\bfa}{\begin{fact} \rm}
\newcommand{\efa}{\end{fact}}
\newcommand{\bp}{\begin{proposition} \rm}
\newcommand{\ep}{\end{proposition}}
\newcommand{\bd}{\begin{definition} \rm}
\newcommand{\ed}{\end{definition}}
\newcommand{\bcs}{\begin{construction} \rm}
\newcommand{\ecs}{\end{construction}}
\newcommand{\bcd}{\begin{condition} \rm}
\newcommand{\ecd}{\end{condition}}
\newcommand{\bt}{\begin{theorem} \rm}
\newcommand{\et}{\end{theorem}}
\newcommand{\bn}{\begin{notation} \rm}
\newcommand{\en}{\end{notation}}
\newcommand{\bfi}{\begin{bild} \rm}
\newcommand{\efi}{\end{bild}}
\newcommand{\bsta}{\begin{statement} \rm}
\newcommand{\esta}{\end{statement}}
\newcommand{\bcom}{\begin{comment} \rm}
\newcommand{\ecom}{\end{comment}}
\newcommand{\bdia}{\begin{diagram} \rm}
\newcommand{\edia}{\end{diagram}}

\newcommand{\bfc}{\begin{figure}[htb] \begin{center}}
\newcommand{\efc}{\end{center} \end{figure}}

\documentstyle[12pt]{article}
\title{
AN ANALYSIS OF DEFEASIBLE INHERITANCE SYSTEMS
}

\author{Karl Schlechta
\thanks{
ks@cmi.univ-mrs.fr, karl.schlechta@web.de, http://www.cmi.univ-mrs.fr/ $\sim$ ks
} \\
Laboratoire d'Informatique Fondamentale de Marseille
\thanks{
UMR 6166, CNRS and Universit\'{e} de Provence,
Address: CMI, 39, rue Joliot-Curie, F-13453 Marseille Cedex 13, France
}
}

\date{August 18, 2007}

\begin{document}

\newtheorem{lemma}{Lemma}[section]
\newtheorem{theorem}[lemma]{Theorem}
\newtheorem{proposition}[lemma]{Proposition}
\newtheorem{corollary}[lemma]{Corollary}
\newtheorem{claim}[lemma]{Claim}
\newtheorem{fact}[lemma]{Fact}
\newtheorem{remark}[lemma]{Remark}
\newtheorem{definition}{Definition}[section]
\newtheorem{construction}{Construction}[section]
\newtheorem{condition}{Condition}[section]
\newtheorem{example}{Example}[section]
\newtheorem{notation}{Notation}[section]
\newtheorem{bild}{Figure}[section]
\newtheorem{comment}{Comment}[section]
\newtheorem{statement}{Statement}[section]
\newtheorem{diagram}{Diagram}[section]

\maketitle

\renewcommand{\labelenumi}
  {(\arabic{enumi})}
\renewcommand{\labelenumii}
  {(\arabic{enumi}.\arabic{enumii})}
\renewcommand{\labelenumiii}
  {(\arabic{enumi}.\arabic{enumii}.\arabic{enumiii})}
\renewcommand{\labelenumiv}
  {(\arabic{enumi}.\arabic{enumii}.\arabic{enumiii}.\arabic{enumiv})}

\begin{abstract}

We give a conceptual analysis of (defeasible or nonmonotonic) inheritance
diagrams, and compare our analysis to the "small/big sets" of preferential
and
related reasoning.

In our analysis, we consider nodes as information sources and truth
values,
direct links as information, and valid paths as information channels and
comparisons of truth values. This results in an upward chaining, split
validity,
off-path preclusion inheritance formalism.

We show that the small/big sets of preferential reasoning have to be
relativized if we want them to conform to inheritance theory, resulting in
a
more cautious approach, perhaps closer to actual human reasoning.

Finally, we interpret inheritance diagrams as theories of prototypical
reasoning, based on two distances: set difference, and information
difference.

We will also see that some of the major distinctions between inheritance
formalisms are consequences of deeper and more general problems of
treating
conflicting information.


\end{abstract}

\tableofcontents

\section{
Introduction
}
\label{Section 1}

Throughout this paper, "inheritance" will stand for "nonmonotonic or
defeasible
inheritance". We will use indiscriminately "inheritance system",
"inheritance
diagram", "inheritance network", "inheritance net".

In this introduction, we first give the motivation for this article, then
describe in very brief terms the basic components of inheritance diagrams,
mention the basic ideas of our analysis, as well as some more general
decisions
about treating contradictory information.

\paragraph{
Motivation
}

$\hspace{0.01em}$

Inheritance sytems or diagrams have an intuitive appeal. They seem close
to
human reasoning, natural, and are also implemented (see  \cite{Mor98}). Yet,
they
are a more procedural approach to nonmonotonic reasoning, and, to the
author's
knowledge, a conceptual analysis, leading to a formal semantics, as well
as a
comparison to more logic based formalisms like the systems $P$ and $R$ of
preferential systems are lacking. We attempt to reduce
the gap between the more procedural and the more analytical approaches in
this
particular case. This will also give indications how to modify systems $P$
and $R$
to approach them more to actual human reasoning. Moreover, we establish a
link
to multi-valued logics and the logics of information sources (see e.g.
 \cite{ABK07} and forthcoming work of the same authors, and also
 \cite{BGH95}).

\paragraph{
Inheritance diagrams:
}

$\hspace{0.01em}$

Inheritance diagrams are deceptively simple. Their conceptually
complicated
nature is seen by e.g. the fundamental difference between direct links and
valid paths, and the multitude of existing formalisms, upward vs. downward
chaining, intersection of extensions vs. direct scepticism, on-path vs.
off-path
preclusion (or pre-emption), split validity vs. total validity preclusion
etc.,
to name a few. Such a proliferation of formalisms usually hints at deeper
problems on the conceptual side, i.e. that the underlying ideas are
ambigous,
and not sufficiently analysed. Therefore, any clarification and resulting
reduction of possible formalisms seems a priori to make progress. Such
clarification will involve conceptual decisions, which need not be shared
by
all, they can only be suggestions. Of course, a proof that such decisions
are
correct is impossible, and so is its contrary.

\paragraph{
Our analysis:
}

$\hspace{0.01em}$

We will introduce into the analysis of inheritance systems a number of
concepts
not usually found in the field, like multiple truth values, access to
information, comparison of truth values, etc. We think that this
additional
conceptual burden pays off by a better comprehension and analysis of the
problems behind the surface of inheritance.

We will also see that some distinctions between inheritance formalisms go
far
beyond questions of inheritance, and concern general problems of treating
contradictory information - isolating some of these is another objective
of
this article.

The text is essentially self-contained, still some familiarity with the
basic
concepts of inheritance systems and nonmonotonic logics in general is
helpful.
For a presentation, the reader might look into  \cite{Sch97-2} and
 \cite{Sch04}.

A.Bochman, has pointed out to author work by J.Barwise, D.Gabbay,
C.Hartonas,  \cite{BGH95}, on information flow. This has a
superficial resemblance
with the present pages. But, first, the BGH work is much deeper into
logic,
presenting sequent calculi, completeness results, etc. Second, our work is
on non-monotonic logics, which BGH is not, and our main thrust is a
conceptual
analysis of inheritance networks, also working with multiple truth values.
But
the basic ideas are about similar situations.

The text is organized as follows. After an introduction to inheritance
theory
and big/small subsets and the systems $P$ and $R$ in Section 2 and Section
3, we
turn to an informal description of the fundamental differences between
inheritance and the systems $P$ and $R$ in Section 4.2, give an analysis
of
inheritance systems in terms of information and information flow in
Section 4.3,
then in terms of reasoning with prototypes in Section 4.4, and conclude in
Section 5 with a translation of inheritance into (necessarily deeply
modified)
coherent systems of big/small sets, respectively logical systems $P$ and
$R.$

One of the main modifications will be to relativize the notions of
small/big,
which thus become less "mathematically pure" but perhaps closer to actual
use in "dirty" common sense reasoning.
\section{
Introduction to nonmonotonic inheritance
}
\label{Section 2}
\subsection{
Basic discussion
}
\label{Section 2.1}

We give here an informal discussion. The reader unfamiliar with
inheritance
systems should consult in parallel Definition \ref{Definition 2.3} and
Definition \ref{Definition 2.4}. As there are many variants of the
definitions, it seems
reasonable to discuss them before a formal introduction, which, otherwise,
would pretend to be definite without being so.

\paragraph{
(Defeasible or nonmonotonic) inheritance networks or diagrams
}

$\hspace{0.01em}$

Nonmonotonic inheritance systems describe situations like
"normally, birds fly", written $birds \xcp fly.$ Exceptions are
permitted, "normally penguins $don' t$ fly", $penguins \xcP fly.$

\bd

$\hspace{0.01em}$

(+++*** Orig. No.:  Definition 2.1 )

\label{Definition 2.1}

A nonmonotonic inheritance net is a finite DAG, directed, acyclic graph,
with two types of arrows or links, $ \xcp $ and $ \xcP,$ and labelled
nodes. We will use
$ \xbG $ etc. for such graphs, and $ \xbs $ etc. for paths - the latter to
be defined below.

\ed

Roughly (and to be made precise and modified below, we try to give here
just a
first intuition), $X \xcp Y$ means that
"normal" elements of $X$ are in $Y,$ and $X \xcP Y$ means that "normal"
elements of $X$ are
not in $Y.$ In a semi-quantitative set interpretation, we will read "most
elements
of $X$ are in $Y'',$ "most elements of $X$ are not in $Y'',$ etc. This
is by no means
the only interpretation, as we will see.

According to the set intrepretation,
we will also use informally expressions like $X \xcs Y,$ X-Y, $ \xdC X,$
etc. But we
will also use nodes informally as formulas, like $X \xcu Y,$ $X \xcu \xCN
Y,$ $ \xCN X,$ etc.
All this will only be used as an appeal to intuition.

Nodes at the beginning of an arrow can also stand for individuals,
so $Tweety \xcP fly$ means something like: "normally, Tweety will not
fly". As always
in nonmonotonic systems, exceptions are permitted, so the soft rules
"birds
fly", "penguins $don' t$ fly", and (the hard rule) "penguins are birds"
can coexist
in one diagram, penguins are then abnormal birds (with respect to flying).
The direct link $penguins \xcP fly$ will thus be accepted, or considered
valid, but
not the composite path $penguins \xcp birds \xcp fly,$ by specificity -
see below.
This is illustrated by Diagram \ref{Diagram 2.1},
where a stands for Tweety, $c$ for penguins,
$b$ for birds, $d$ for flying animals or objects.

(Remark: The arrows $a \xcp c,$ $a \xcp b,$ and $c \xcp b$ can also be
composite
paths - see below for the details.)

\bfc

\bdia

\label{Diagram 2.1}

\unitlength1.0mm
\begin{picture}(130,100)

\newsavebox{\SECHSacht}
\savebox{\SECHSacht}(140,110)[bl]
{

\put(0,95){The Tweety diagram}

\put(43,27){\vector(1,1){24}}
\put(37,27){\vector(-1,1){24}}
\put(13,57){\vector(1,1){24}}
\put(67,57){\vector(-1,1){24}}

\put(53,67){\line(1,1){4}}

\put(67,54){\vector(-1,0){54}}

\put(39,23){$a$}
\put(9,53){$b$}
\put(69,53){$b$}
\put(39,83){$d$}

}

\put(0,0){\usebox{\SECHSacht}}
\end{picture}

\edia

\efc

(Of course, there is an analogous case for the opposite polarity, i.e.
when the arrow from $b$ to $d$ is negative, and the one from $c$ to $d$ is
positive.)

The main problem is to define in an intuitively acceptable way
a notion of valid path, i.e. concatenations of arrows satisfying certain
properties.

We will write $ \xbG \xcm \xbs,$ if $ \xbs $ is a valid path in the
network $ \xbG,$ and if
$x$ is the origin, and $y$ the endpoint of $ \xbs,$ and $ \xbs $ is
positive, we will write
$ \xbG \xcm xy,$ i.e. we will accept the conclusion that x's are $y' $s,
and analogously
$ \xbG \xcm x \ol{y}$ for negative paths. Note that we
will not accept any other conclusions, only those established by a valid
path, so many questions about conclusions have a trivial negative answer:
there
is obviously no path from $x$ to $y.$ E.g., there is no path from $b$ to
$c$ in
Diagram \ref{Diagram 2.1}. Likewise, there are no disjunctions,
conjunctions etc. in our conclusions, and negation is present only in a
strong
form: "it is not the case that x's are normally $y' $s" is not a possible
conclusion, only "x's are normally not $y' $s" is one. Also, possible
contradictions are contained, there is no EFQ.

To simplify matters, we assume that for no two nodes $x,y \xbe \xbG $ $x
\xcp y$ and $x \xcP y$ are
both in $ \xbG,$ intuitively, that $ \xbG $ is free from (hard)
contradictions. This
restriction is inessential for our purposes. We admit, however, soft
contradictions, and preclusion, which allows us to solve some soft
contradictions -
as we already did in the penguins example. We will also assume that all
arrows
stand for rules with possibly exceptions, again, this restriction is not
important for our purposes. Moreover, in the abstract treatment, we will
assume that all nodes stand for sets, though this will not be true for all
examples discussed.

This might be the place for a remark on absence of cycles. Suppose we also
have a positive arrow from $b$ to $c$ in Diagram \ref{Diagram 2.1}.
Then, the concept of
preclusion collapses,
as there are now equivalent arguments to accept $a \xcp b \xcp d$ and $a
\xcp c \xcP d.$ Thus,
if we do not want to introduce new complications, we cannot rely on
preclusion
to decide conflicts. It seems that this will change the whole outlook on
such diagrams.

Inheritance networks were introduced about 20 years ago (see e.g.  \cite{Tou84},
 \cite{Tou86},  \cite{THT87}), and exist
in a multitude of more or less differing formalisms, see
e.g.  \cite{Sch97-2} for a brief discussion. There still does not
seem
to exist a satisfying semantics for these networks. The author's own
attempt
 \cite{Sch90} is an a posteriori semantics, which cannot explain or
criticise or
decide between the different formalisms. We will give here a conceptual
analysis, which provides also at least some building blocks for a
semantics,
and a translation into (a modified version of) the language of small/big
subsets, familiar from preferential structures, see below in Section 3.

We will now discuss the two fundamental situations of contradictions,
then give a detailed inductive definition of valid paths for a certain
formalism so the reader has firm ground under his feet, and then present
briefly
some alternative formalisms.

As in all of nonmonotonic reasoning, the interesting questions arise in
the
treatment of contradictions and exceptions. The
difference in quality of information is expressed by "preclusion"
(or $'' pre-emption'' ).$ The basic diagram is the Tweety diagram,
see Diagram \ref{Diagram 2.1}.

Unresolved contradictions give either rise to a branching
into different extensions, which may roughly be seen as maximal consistent
subsets, or to mutual cancellation in directly sceptical approaches.
The basic diagram for the latter is the Nixon Diamond, see
Diagram \ref{Diagram 2.2}, where $a=Nixon,$ $b=Quaker,$ $c=Republican,$
$d=pacifist.$

In the directly sceptical approach, we will not accept any path from a to
$d$
as valid, as there is an unresolvable contradiction between the two
candidates.

\bfc

\bdia

\label{Diagram 2.2}

\unitlength1.0mm
\begin{picture}(130,100)

\newsavebox{\ZWEIzwei}
\savebox{\ZWEIzwei}(140,110)[bl]
{

\put(0,95){The Nixon Diamond}

\put(43,27){\vector(1,1){24}}
\put(37,27){\vector(-1,1){24}}
\put(13,57){\vector(1,1){24}}
\put(67,57){\vector(-1,1){24}}

\put(53,67){\line(1,1){4}}

\put(39,23){$a$}
\put(9,53){$b$}
\put(69,53){$c$}
\put(39,83){$d$}

}

\put(0,0){\usebox{\ZWEIzwei}}
\end{picture}

\edia

\efc

The extensions approach can be turned into an indirectly sceptical one,
by forming first all extensions, and then taking the intersection of
either the
sets of valid paths, or of valid conclusions, see [MS91] for a detailed
discussion.

In more detail:

\paragraph{
Preclusion
}

$\hspace{0.01em}$

In the above example, our intuition tells us that it is not admissible
to conclude from the fact that penguins are birds, and that most birds
fly that most penguins fly. The horizontal arrow $c \xcp b$ together with
$c \xcP d$ barrs
this conclusion, it expresses specificity. Consequently, we have
to define the conditions under which two potential paths neutralize each
other, and when one is victorious. The idea is as follows:
1) We want to be sceptical, in the sense that we do not believe every
potential path. We will not arbitrarily chose one either.
2) Our scepticism will be restricted, in the sense that we will often make
well defined choices for one path in the case of conflict:
a) If a compound potential path is in conflict with a direct link,
the direct link wins.
$b)$ Two conflicting paths of the same type neutralize each other, as in
the
Nixon Diamond, where neither potential path will be valid.
$c)$ More specific information will win over less specific one.

(It is essential in the Tweety diagram that the arrow $c \xcP d$ is a
direct link, so
it is in a
way stronger than compound paths.) The arrows $a \xcp b,$ $a \xcp c,$ $c
\xcp b$ can also
be composite paths: The path from $c$ to $b$ (read $c \xcc $  \Xl  $ \xcc
b$!),
however, tells us, that the information coming from $c$ is more specific
(and thus considered more reliable), so the negative path from a to $d$
via
$c$ will win over the positive one via $b.$ The precise inductive
definition
will be given below. This concept is evidently independent of the lenght
of the
paths, $a \xFB \xcp c$ may be much longer than $a \xFB \xcp b,$ so this is
not shortest
path reasoning (which has some nasty drawbacks, discussed e.g. in
[HTT87]).

Before we give a formalism based on these ideas, we refine them, adopt one
possibility (but indicate some modifications), and discuss alternatives
later.
\subsection{
Directly sceptical split validity upward chaining off-path inheritance
}
\label{Section 2.2}

Our approach will be directly sceptical, i.e. unsolvable contradictions
result
in the absence of valid paths, it is upward chaining, and split-validity
for
preclusions (discussed below). We will indicate modifications to make it
extension based, as well as for total-validity preclusion.
This approach is strongly inspired by classical work in the field
by Horty, Thomason, Touretzky, and others, and we claim no priority
whatever.
If it is new at all, it is a very minor modification of existing
formalisms.

Our conceptual ideas to be presented in detail in Section 4.3 make split
validity, off-path preclusion and upward chaining a natural choice.

For the reader's convenience, we give here a very short resume of these
ideas:
We consider only arrows as information, e.g. $a \xcp b$ will be considered
information $b$ valid at or for a. Valid (composed positive) paths will
not be
considered (direct) information. They will be seen as a way to obtain
information, so a valid path $ \xbs:x \Xl  \xcp a$ makes information $b$
accessible to $x,$
and, secondly, as a means of comparing information
strength, so a valid path $ \xbs:a \Xl. \xcp a' $ will make information
at a stronger than
information at $a'.$ Valid negative paths have no function, we will only
consider
the positive initial part as discussed above, and the negative end arrow
as
information, but never the whole path.

Choosing direct scepticism is a decision beyond the scope of this article,
and
we just make it. It is a general question how to treat contradictory and
absent
information, and if they are equivalent or not, see the remark in
Section \ref{Section 4.4}. (The fundamental difference
between intersection of extensions and direct scepticism for defeasible
inheritance was shown in  \cite{Sch93}.)

We turn now to the announced variants as well as a finer distinction
within the
directly sceptical approach.

Our approach generates another problem, essentially that of the treatment
of a
mixture of contradictory and concordant information of multiple strengths
or
truth values. We bundle the decision of this problem with that for direct
scepticism into a "plug-in" decision, which will be used in three
approaches:
the conceptual ideas, the inheritance algorithm, and the choice of the
reference
class for subset size (and implicitly also for the treatment as a
prototype
theory). It is thus well encapsulated, and independent from the context.

These decisions (but, perhaps to a lesser degree, (1)) concern a wider
subject
than only inheritance
networks. Thus, it is not surprising that there are different formalisms
for
solving such networks, deciding one way or the other. But this multitude
is not
the fault of inheritance theory, it is only a symptom of a deeper
question. We
first give an overview for a clearer overall picture,
and discuss them in detail below, as they involve sometimes quite subtle
questions.

(1) Upward chaining against downward or double chaining.

(2.1) Off-path against on-path preclusion.

(2.2) Split validity preclusion against total validity preclusion.

(3) Direct scepticism against intersection of extensions.

(4) Treatment of mixed contradiction and preclusion situations,
no preclusion by paths of the same polarity.

(1) This can also be seen as a difference in reasoning from cause to
effect vs.
backward reasoning, looking for causes for an effect. (A word of warning:
There
is a
well-known article  \cite{SL89} from which a superficial reader
might conclude that
upward chaining is tractable, and downward chaining is not. A more careful
reading reveals that, on the negative side, the authors only show that
double
chaining is not tractable.) We will adopt upward chaining in all our
approaches.
See Section \ref{Section 4.4} for more remarks.

(2.1) and (2.2) Both are consequences of our view - to be discussed below
in
Section 4.3 - to see valid paths also as
an absolute comparison of truth values, independent of reachability of
information. This question of absoluteness transcends obviously
inheritance
networks. Our decision is, of course, again uniform for all our
approaches.

(3) This point, too, is much more general than the problems of
inheritance.
It is, among other things, a question of whether only the two possible
cases
(positive and negative) may hold, or whether there might be still other
possibilities. See Section \ref{Section 4.4}.

(4) This concerns the treatment of truth values in more complicated
situations, where we have a mixture of agreeing and contradictory
information.
Again, this problem reaches far beyond inheritance networks.

We will group (3) and (4) together in one general, "plug-in"
decision, to be found in all approaches we discuss.

\bd

$\hspace{0.01em}$

(+++*** Orig. No.:  Definition 2.2 )

\label{Definition 2.2}

A plug-in decision:

We describe now more precisely a situation which we will meet in all
contexts
discussed, and
whose decision goes beyond our problem - thus, we have to adopt one or
several
alternatives, and translate them into the approaches we will discuss.
There will be one global decision, which is (and can be) adapted to the
different contexts.

Suppose we have information about $ \xbf $ and $ \xbq,$ where $ \xbf $
and $ \xbq $ are presumed to be
independent - in some adequate sense.

Suppose then that we have information sources $A_{i}:i \xbe I$ and
$B_{j}:j \xbe J,$ where
the $A_{i}$ speak about $ \xbf $ (they say $ \xbf $ or $ \xCN \xbf ),$ and
the $B_{j}$ speak about $ \xbq $ in the
same way. Suppose further that we have a partial, not necessarily
transitive (!),
ordering $<$ on the information sources $A_{i}$ and $B_{j}$ together.
$X<Y$ will say that $X$ is
better (intuition: more specific) than $Y.$ (The potential lack of
transitivity
is crucial, as valid paths do not always concatenate to valid paths - just
consider the Tweety diagram.)

We also assume that there are contradictions, i.e. some $A_{i}$ say $ \xbf
,$ some $ \xCN \xbf,$
likewise for the $B_{j}$ - otherwise, there are no problems in our
context.

We can now take several approaches, all taking contradictions and the
order $<$
into account.
 \xEI
 \xDH (P1) We use the global relation $<,$ and throw away all information
coming from
sources of minor quality, i.e. if there is $X$ such that $X<Y,$ then no
information
coming from $Y$ will be taken into account. Consequently, if $Y$ is the
only source
of information about $ \xbf,$ then we will have no information about $
\xbf.$ This seems
an overly radical approach, as one source might be better for $ \xbf,$
but not
necessarily for $ \xbq,$ too.

If we adopt this approach, we can continue as below, and can even split
in analogue ways into (P1.1) and (P1.2), as we do below for (P2.1) and
(P2.2).

 \xDH (P2) We consider the information about $ \xbf $ separately from the
information about
$ \xbq.$ Thus, we consider for $ \xbf $ only the $A_{i},$ for $ \xbq $
only the $B_{j}.$ Take now e.g. $ \xbf $
and the $A_{i}.$ Again, there are (at least) two alternatives.
 \xEI
 \xDH (P2.1) We eliminate again all sources among the $A_{i}$ for which
there is a better
$A_{i' },$ irrespective of whether they agree on $ \xbf $ or not.
 \xEI
 \xDH (a) If the sources left are contradictory, we conclude nothing about
$ \xbf,$ and
accept for $ \xbf $ none of the sources. (This is a directly sceptical
approach of
treating unsolvable contradictions, following our general strategy.)

 \xDH (b) If the sources left agree for $ \xbf,$ i.e. all say $ \xbf,$
or all say $ \xCN \xbf,$ then we
conclude $ \xbf $ (or $ \xCN \xbf ),$ and accept for $ \xbf $ all the
remaining sources.
 \xEJ
 \xDH (P2.2) We eliminate again all sources among the $A_{i}$ for which
there is a better
$A_{i' }$, but only if $A_{i}$ and $A_{i' }$ have contradictory
information.
Thus, more sources may survive than in approach (P2.1).

We now continue as for (P2.1):
 \xEI
 \xDH (a) If the sources left are contradictory, we conclude nothing about
$ \xbf,$ and
accept for $ \xbf $ none of the sources.

 \xDH (b) If the sources left agree for $ \xbf,$ i.e. all say $ \xbf,$
or all say $ \xCN \xbf,$ then we
conclude $ \xbf $ (or $ \xCN \xbf ),$ and accept for $ \xbf $ all the
remaining sources.
 \xEJ
 \xEJ
 \xEJ
The difference between (P2.1) and (P2.2) is illustrated by the following
simple
example. Let $A<A' <A'',$ but $A \xEc A'' $ (recall that $<$ is not
necessarily transitive),
and $A \xcm \xbf,$ $A' \xcm \xCN \xbf,$ $A'' \xcm \xCN \xbf.$ Then
(P2.1) decides for $ \xbf $ $( \xCB $ is the only
survivor), (P2.2) does not decide, as $ \xCB $ and $ \xCB '' $ are
contradictory, and both
survive in (P2.2).

There are arguments for and against either solution: (P2.1) gives a
uniform
picture, more independent from $ \xbf,$ (P2.2) gives more weight to
independent
sources, it "adds" information sources, and thus gives potentially more
weight
to information from several sources. (P2.2) seems more in the tradition of
inheritance networks, so we will consider it in the further development.

The reader should note that our approach is quite far from a fixed point
approach in two ways: First, fixed point approaches seem more appropriate
for extensions based approaches, as both try to collect a maximal set of
uncontradictory information. Second, we eliminate information when there
is better, contradicting information, even if the final result agrees with
the first. This, too, contradicts in spirit the fixed point approach.

\ed

After these preparations, we turn to a formal definition of validity of
paths.

\paragraph{
The definition of $ \xcm $ (i.e. of validity of paths)
}

$\hspace{0.01em}$

All definitions are relative to a fixed diagram $ \xbG.$
The notion of degree will be defined relative to all nodes of $ \xbG,$ as
we will
work with split validity preclusion, so the paths to consider may have
different
origins. For simplicity, we consider $ \xbG $ to be just a set of points
and
arrows, thus e.g. $x \xcp y \xbe \xbG $ and $x \xbe \xbG $ are defined,
when $x$ is a point in $ \xbG,$ and
$x \xcp y$ an arrow in $ \xbG.$ Recall that we have two types of arrows,
positive and
negative ones.

We first define generalized and potential paths, then the notion of
degree, and
finally validity of paths, written $ \xbG \xcm \xbs,$ if $ \xbs $ is a
path, as well as $ \xbG \xcm xy,$
if $ \xbG \xcm \xbs $ and $ \xbs:x \Xl. \xcp y.$

\bd

$\hspace{0.01em}$

(+++*** Orig. No.:  Definition 2.3 )

\label{Definition 2.3}

(1) Generalized paths:

A generalized path is an uninterrupted chain of positive or negative
arrows
pointing in the same direction, more precisely:

$x \xcp p \xbe \xbG $ $ \xcp $ $x \xcp p$ is a generalized path,

$x \xcP p \xbe \xbG $ $ \xcp $ $x \xcP p$ is a generalized path.

If $x \xFB \xcp p$ is a generalized path, and $p \xcp q \xbe \xbG $,
then
$x \xFB \xcp p \xcp q$ is a generalized path,

if $x \xFB \xcp p$ is a generalized path, and $p \xcP q \xbe \xbG $,
then
$x \xFB \xcp p \xcP q$ is a generalized path.

(2) Concatenation:

If $ \xbs $ and $ \xbt $ are two generalized paths, and the end point of $
\xbs $ is the same
as the starting point of $ \xbt,$ then $ \xbs \xDM \xbt $ is the
concatenation of $ \xbs $ and $ \xbt.$

(3) Potential paths (pp.):

A generalized path, which contains at most one negative arrow, and then at
the
end, is a potential path. If the last link is positive, it is a positive
potential path, if not, a negative one.

(4) Degree:

As already indicated, we shall define paths inductively. As we do not
admit cycles in our systems, the arrows define a well-founded relation
on the vertices. Instead of using this relation for
the induction, we shall first define the auxiliary notion of
degree, and do induction on the degree.
Given a node $x$ (the origin), we need a (partial) mapping $f$ from the
vertices to
natural numbers such that $p \xcp q$ or $p \xcP q$ $ \xbe $ $ \xbG $
implies $f(p)<f(q),$ and define (relative
to $x):$

Let $ \xbs $ be a generalized path from $x$ to $y,$ then $deg_{ \xbG,x}(
\xbs ):=deg_{ \xbG,x}(y):=$
the maximal length of any generalized path parallel to $ \xbs,$ i.e.
beginning in $x$ and ending in $y.$

\ed

\bd

$\hspace{0.01em}$

(+++*** Orig. No.:  Definition 2.4 )

\label{Definition 2.4}

Inductive definition of $ \xbG \xcm \xbs:$

Let $ \xbs $ be a potential path.
 \xEI
 \xDH Case $I$:

$ \xbs $ is a direct link in $ \xbG.$ Then $ \xbG \xcm \xbs $

(Recall that we have no hard contradictions in $ \xbG.)$

 \xDH Case II:

$ \xbs $ is a compound potential path, $deg_{ \xbG,a}( \xbs )=n,$ and $
\xbG \xcm \xbt $ is defined
for all $ \xbt $ with degree less than $n$ - whatever their origin and
endpoint.

 \xDH Case II.1:

Let $ \xbs $ be a positive pp. $x \xFB \xcp u \xcp y,$ let $ \xbs ':=x
\xFB \xcp u,$ so $ \xbs = \xbs ' \xDM u \xcp y$

Then, intuitively, $ \xbG \xcm \xbs $ iff
 \xEh
 \xDH (1) $ \xbs $ is a candidate by upward chaining,
 \xDH (2) $ \xbs $ is not precluded by more specific contradicting
information,
 \xDH (3) all potential contradictions are themselves precluded by
information
contradicting them.
 \xEj
Note that (2) and (3) are the translation of (P2.2) in Definition
\ref{Definition 2.2}.

Formally, $ \xbG \xcm \xbs $ iff
 \xEh
 \xDH (1) $ \xbG \xcm \xbs ' $ and $u \xcp y \xbe \xbG.$

(The initial segment must be a path, as we have an upward chaining
approach.
This is decided by the induction hypothesis.)

 \xDH (2) There are no $v,$ $ \xbt,$ $ \xbt ' $ such that $v \xcP y \xbe
\xbG $ and $ \xbG \xcm \xbt:=x \xFB \xcp v$ and
$ \xbG \xcm \xbt ':=v \xFB \xcp u.$ $( \xbt $ may be the empty path, i.e.
$x=v.)$

$( \xbs $ itself is not precluded by split validity preclusion and a
contradictory
link. Note that $ \xbt \xDM v \xcP y$ need not be valid, it suffices
that it is a better candidate (by $ \xbt ' ).)$

 \xDH (3) all potentially conflicting paths are precluded by information
contradicting them:

For all $v$ and $ \xbt $ such that $v \xcP y \xbe \xbG $ and $ \xbG \xcm
\xbt:=x \xFB \xcp v$ (i.e. for all potentially
conflicting paths $ \xbt \xDM v \xcP y)$ there is $z$ such that $z \xcp y
\xbe \xbG $ and either

$z=x$

(the potentially conflicting pp. is itself precluded by a direct link,
which is
thus valid)

or

there are $ \xbG \xcm \xbr:=x \xFB \xcp z$ and $ \xbG \xcm \xbr ':=z
\xFB \xcp v$ for suitable $ \xbr $ and $ \xbr '.$
 \xEj
 \xDH Case II.2: The negative case, i.e.
$ \xbs $ a negative pp. $x \xFB \xcp u \xcP y,$ $ \xbs ':=x \xFB \xcp u,$
$ \xbs = \xbs ' \xDM u \xcP y$
is entirely symmetrical.
 \xEJ

\ed

\br

$\hspace{0.01em}$

(+++*** Orig. No.:  Remark 2.1 )

\label{Remark 2.1}

The following remarks all concern preclusion.

(1) Thus, in the case of preclusion, there is a valid path from $x$ to
$z,$
and $z$ is more specific than $v,$ so
$ \xbt \xDM v \xcP y$ is precluded. Again, $ \xbr \xDM z \xcp y$ need not
be a valid path, but it is
a better candidate than $ \xbt \xDM v \xcP y$ is, and as $ \xbt \xDM v
\xcP y$ is in simple contradiction,
this suffices.

(2) Our definition is stricter than many usual ones, in the
following sense: We require - according to our general picture to treat
only
direct links as information - that the preclusion "hits" the precluded
path
at the end, i.e. $v \xcP y \xbe \xbG,$ and $ \xbr ' $ hits $ \xbt \xDM v
\xcP y$ at $v.$ In other definitions,
it is possible that the preclusion hits at some $v',$ which is somewhere
on the
path $ \xbt,$ and not necessarily at its end. For instance, in the Tweety
Diagram,
see Diagram \ref{Diagram 2.1}, if there were a node $b' $ between $b$
and $d,$ we will need the
path $c \xcp b \xcp b' $ to be valid, (obvious) validity of the arrow $c
\xcp b$ will not
suffice.

(3) If we allow $ \xbr $ to be the empty path, then the case $z=x$ is a
subcase of the
present one.

(4) Our conceptual analysis has led to a very important simplification of
the
definition of validity. If we adopt on-path preclusion, we have to
remember
all paths which led to the information source to be considered: In the
Tweety
diagram, we have to remember that there is an arrow $a \xcp b,$ it is not
sufficient
to note that we somehow came from a to $b$ by a valid path, as the path $a
\xcp c \xcp b \xcp d$
is precluded, but not the path $a \xcp b \xcp d.$ If we adopt total path
preclusion, we
have to remember the valid path $a \xcp c \xcp b$ to see that it precludes
$a \xcp c \xcp d.$ If
we allow preclusion to "hit" below the last node, we also have to remember
the entire path which is precluded. Thus, in all those cases, whole paths
(which can be very long) have to be remembered, but NOT in our definition.

We only need to remember (consider the Tweety diagram):
(a) we want to know if $a \xcp b \xcp d$ is valid, so we have to remember
a, $b,$ $d.$
Note that the (valid) path from a to $b$ can be composed and very long.
(b) we look at possible preclusions, so we have to remember $a \xcp c \xcP
d,$ again
the (valid) path from a to $c$ can be very long.
(c) we have to remember that the path from $c$ to $b$ is valid (this was
decided
by induction before).

So in all cases (the last one is even simpler), we need only remember the
starting node, a (or $c),$ the last node of the valid paths, $b$ (or $c),$
and the
information $ \xcp d$ or $ \xcP d$ - i.e. the size of what has to be
recalled is
$ \xck 3.$ (Of course, there may be many possible preclusions, but in all
cases we have to look at a very limited situation, and not arbitrarily
long
paths.)

We take a fast look forward to Section 4.3, where we describe diagrams as
information and its transfer, and nodes also as truth values. In these
terms -
and the reader is asked to excuse the digression - we may note above point
(a) as $a \xch_{b}d$ - expressing that, seen from a, $d$ holds with truth
value $b,$
(b) as $a \xch_{c} \xCN d,$ (c) as $c \xch_{c}b$ - and this is all we need
to know.

$ \xcz $
\\[3ex]

\er

We indicate here some modifications of the definition without discussion,
which
is to be found below.

(1) For on-path preclusion only:
Modify condition (2) in Case II.1 to:
$(2' )$ There is no $v$ on the path $ \xbs $ (i.e. $ \xbs:x \xFB \xcp v
\xFB \xcp u)$ such that $v \xcP y \xbe \xbG.$

(2) For total validity preclusion:
Modify condition (2) in Case II.1 to:
$(2' )$ There are no $v,$ $ \xbt,$ $ \xbt ' $ such that $v \xcP y \xbe
\xbG $ and $ \xbt:=x \xFB \xcp v$ and
$ \xbt ':=v \xFB \xcp u$ such that $ \xbG \xcm \xbt \xDM \xbt '.$

(3) For extension based approaches:
Modify condition (3) in Case II.1 as follows:
$(3' )$ If there are conflicting paths,
which are not precluded themselves by contradictory information, then we
branch recursively (i.e. for all such situations) into two extensions,
one,
where the positive non-precluded paths are valid, one, where the negative
non-precluded paths are valid.

\bd

$\hspace{0.01em}$

(+++*** Orig. No.:  Definition 2.5 )

\label{Definition 2.5}

Finally, define $ \xbG \xcm xy$ iff there is $ \xbs:x \xcp y$ s.th. $
\xbG \xcm \xbs,$ likewise for $x \ol{y}$ and
$ \xbs:x \xFB \xcP y.$

\ed

Diagram \ref{Diagram 2.3} shows the most complicated situation for the
positive case.

\bfc

\bdia

\label{Diagram 2.3}

\unitlength1.0mm
\begin{picture}(130,100)

\newsavebox{\Preclusion}
\savebox{\Preclusion}(140,90)[bl]
{


\multiput(43,8)(1,1){5}{\circle*{.3}}
\put(48,13){\vector(1,1){17}}

\multiput(37,8)(-1,1){5}{\circle*{.3}}
\put(32,13){\vector(-1,1){17}}

\put(13,38){\vector(1,1){24}}
\put(67,38){\vector(-1,1){24}}

\put(40,37){\vector(0,1){23}}

\multiput(40,8)(0,1){5}{\circle*{.3}}
\put(40,13){\vector(0,1){17}}

\multiput(66,34)(-1,0){5}{\circle*{.3}}
\put(61,34){\vector(-1,0){17}}

\put(39,3){$x$}
\put(9,33){$u$}
\put(39,33){$v$}
\put(69,33){$z$}
\put(39,63){$y$}

\put(38,50){\line(1,0){3.7}}

}

\put(0,0){\usebox{\Preclusion}}
\end{picture}

\edia

\efc

We have to show now that the above approach corresponds to the preceeding
discussion.

\bfa

$\hspace{0.01em}$

(+++*** Orig. No.:  Fact 2.2 )

\label{Fact 2.2}

The above definition and the one outlined in Definition \ref{Definition 2.2}
correspond.

\efa

\paragraph{
Proof:
}

$\hspace{0.01em}$

As Definition \ref{Definition 2.2} is informal, this cannot be a formal
proof, but it
is obvious how to transform it into one.

We argue for the result, the argument for valid paths is similar.

Consider then case (P2.2) in Definition \ref{Definition 2.2},
and start from some $x.$

\paragraph{
Case 1:
}

$\hspace{0.01em}$

Direct links, $x \xcp z$ or $x \xcP z.$

By definition, as a direct link starts at $x,$ the information $z$ or $
\xCN z$ is
stronger than all other accessible information. Thus, the link and the
information will be valid in both approaches. Note that we assumed $ \xbG
$ free from
hard contradictions.

\paragraph{
Case 2:
}

$\hspace{0.01em}$

Composite paths.

In both approaches, the initial segment has to be valid, as information
will
otherwise not be accessible. Also, in both approaches, information will
have
the form of direct links from the accessible source. Thus, condition (1)
in
Case II.1 corresponds to condition (1) in Definition \ref{Definition 2.2}.

In both approaches, information contradicted by a stronger source
(preclusion)
is discarded, as well as information which is contradicted by other, not
precluded sources, so (P2.2) in Definition \ref{Definition 2.2} and II.1
$(2)+(3)$ correspond. Note that variant (P2.1) of Definition \ref{Definition
2.2} would give a
different result - which we could, of course, also imitate in a modified
inheritance approach.

\paragraph{
Case 3:
}

$\hspace{0.01em}$

Other information.

Inheritance nets give no other information, and, as pointed out,
we do not add any other information either in the approach in
Definition \ref{Definition 2.2}.

Thus, both approaches are equivalent.

$ \xcz $
\\[3ex]
\subsection{
Review of other approaches and problems
}
\label{Section 2.3}

We now discuss shortly in more detail some of the differences between
various
major definitions of inheritance formalisms.

Diagram 6.8, $p.$ 179, in  \cite{Sch97-2} (which is probably due to
folklore of the
field) shows requiring downward chaining would be wrong. We repeat it
here,
see Diagram \ref{Diagram 2.4}.

\bfc

\bdia

\label{Diagram 2.4}

\unitlength1.0mm
\begin{picture}(130,100)

\newsavebox{\ZWEIvier}
\savebox{\ZWEIvier}(140,110)[bl]
{

\put(0,95){The problem of downward chaining:}

\put(43,27){\vector(1,1){24}}
\put(37,27){\vector(-1,1){24}}
\put(13,57){\vector(1,1){24}}
\put(67,57){\vector(-1,1){24}}

\put(53,67){\line(1,1){4}}

\put(67,54){\vector(-1,0){54}}

\put(40,7){\vector(0,1){14}}
\put(43,7){\line(3,5){24}}
\put(58,28.1){\line(-5,3){3.6}}

\put(39,3){$z$}
\put(39,23){$u$}
\put(9,53){$v$}
\put(69,53){$x$}
\put(39,83){$y$}

}

\put(0,0){\usebox{\ZWEIvier}}
\end{picture}

\edia

\efc

Preclusions valid above (here at $u)$ can be invalid at lower points (here
at $z),$
as part of the relevant information is not any more accessible (or becomes
accessible). We have $u \xcp x \xcP y$ valid, by downward chaining, any
valid path
$z \xcp u \Xl.y$ has to have a valid final segment u \Xl y, which can
only be $u \xcp x \xcP y,$
but intuition says that $z \xcp u \xcp v \xcp y$ should be valid. Downward
chaining prevents such changes, and thus seems inadequate, so we decide
for
upward chaining. (Already preclusion itself underlines upward chaining: In
the
Tweety diagram, we have to know that the path from bottom up to penguins
is
valid. So at least some initial subpaths have to be known - we need upward
chaining.) (The rejection of downward chaining seems at first sight to be
contrary to the intuitions carried by the word $'' inheritance''.)$
See the remark in Section \ref{Section 4.4}.

\paragraph{
Extension-based versus directly skeptical definitions
}

$\hspace{0.01em}$

As this distinction has already received detailed discussion in the
literature, we shall be very brief here.
An extension of a net is essentially a maximally consistent and in some
appropriate sense reasonable subset of all its potential paths. This can
of course be presented either as a liberal conception (focussing on
individual extensions) or as a skeptical one (focussing on their
intersection - or, the
intersection of their conclusion sets). The seminal presentation is that
of [Tou86], as refined by [San86]. The directly skeptical approach seeks
to
obtain a notion of skeptically accepted path and conclusion, but without
detouring through extensions. Its classic presentation is that of [HTT87].
Even while still searching for fully adequate definitions of either kind,
we may use the former approach as a useful "control" on the latter. For
if we can find an intuitively possible and reasonable extension supporting
a conclusion $x \ol{y},$ whilst a proposed definition for a directly
skeptical
notion of legitimate inference yields xy as a conclusion, then the
counterexemplary extension seems to call into question the adequacy
of the directly skeptical construction, more readily than inversely.

It has been shown in  \cite{Sch93} that the intersection of
extensions is
fundamentally different from the directly sceptical approach.
See the remark in Section \ref{Section 4.4}.

From now on, all definitions considered shall be (at least) upward
chaining.

\paragraph{
On-path versus off-path preclusion
}

$\hspace{0.01em}$

This is a rather technical distinction, discussed in [THT87]. Briefly, a
path
$ \xbs $: $x \xcp  \Xl  \xcp y \xcp  \Xl  \xcp z$ and a direct link $y
\xcP u$ is an off-path preclusion of
$ \xbt $: $x \xcp  \Xl  \xcp z \xcp  \Xl  \xcp u$, but an on-path
preclusion only iff all nodes of $ \xbt $
between $x$ and $z$ lie on the path $ \xbs.$

For instance, in the Tweety diagram, the arrow $c \xcP d$ is an on-path
preclusion
of the path $a \xcp c \xcp b \xcp d,$ but the paths $a \xcp c$ and $c \xcp
b,$ together with $c \xcP d,$ is
an (split validity) off-path preclusion of the path $a \xcp b \xcp c.$

\paragraph{
Split-validity versus total-validity preclusion
}

$\hspace{0.01em}$

Consider again a preclusion $ \xbs:u \xcp  \Xl  \xcp x \xcp  \Xl  \xcp
v,$ and $x \xcP y$ of
$ \xbt:u \xcp  \Xl  \xcp v \xcp  \Xl  \xcp y.$ Most definitions demand
for the preclusion to be
effective - i.e. to prevent $ \xbt $ from being accepted - that the
total path $ \xbs $ is valid. Some ([GV89], [KK89], [KKW89a], [KKW89b])
content
themselves with the combinatorially simpler separate (split) validity of
the lower and upper parts of $ \xbs $: $ \xbs ':u \xcp  \Xl  \xcp x$ and
$ \xbs '':x \xcp  \Xl  \xcp v.$
In Diagram \ref{Diagram 2.5}, taken
from  \cite{Sch97-2}, the path $x \xcp w \xcp v$ is valid,
so is $u \xcp x,$ but not the whole preclusion path $u \xcp x \xcp w \xcp
v.$

\bfc

\bdia

\label{Diagram 2.5}

\unitlength1.0mm
\begin{picture}(130,100)

\newsavebox{\SECHSneun}
\savebox{\SECHSneun}(140,90)[bl]
{

\put(0,75){Split vs. total validity preclusion:}

\put(43,8){\vector(1,1){24}}
\put(37,8){\vector(-1,1){24}}
\put(13,38){\vector(1,1){24}}
\put(67,38){\vector(-1,1){24}}

\put(40,8){\vector(0,1){23}}
\put(66,34){\vector(-1,0){23}}
\put(36,34){\vector(-1,0){23}}

\put(39,3){$u$}
\put(9,33){$v$}
\put(39,33){$w$}
\put(69,33){$x$}
\put(39,63){$y$}

\put(38,20){\line(1,0){3.7}}
\put(52,49){\line(1,1){4}}

}

\put(0,0){\usebox{\SECHSneun}}
\end{picture}

\edia

\efc

Thus, split validity
preclusion will give here the definite result $u \ol{y}.$ With total
validity
preclusion, the diagram has essentially the form of a Nixon Diamond.
\section{
Introduction to small and big sets and the logical systems P and R
}
\label{Section 3}

It is natural to interpret "normality" by some sort of "size": "normality"
might just mean "majority" (perhaps with different weight given to
different
cases), or something like "a big subset". The standard abstraction of
"big"
is the notion of a filter (or, dually, an ideal is the abstraction of
$'' small'' ).$ We include immediately a modification, the weak versions,
to be
discussed below. They seem to be minimal in the following sense: A
reasonable
abstract notion of size without the properties of weak filters seems
difficult
to imagine: The full set seems the best candidate for a "big" subset,
"big"
should cooperate with inclusion, and, finally, no set should be big and
small
at the same time.

\bd

$\hspace{0.01em}$

(+++*** Orig. No.:  Definition 3.1 )

\label{Definition 3.1}

Fix a base set $X.$

A (weak) filter on or over $X$ is a set $ \xdf \xcc \xdp (X)$ - $ \xdp
(X)$ the power set of $X$ -,
such that $(F1)-(F3)$ $((F1),$ (F2), $(F3' )$ respectively) hold:

(F1) $X \xbe \xdf $

(F2) $A \xcc B \xcc X,$ $A \xbe \xdf $ imply $B \xbe \xdf $

(F3) $A,B \xbe \xdf $ imply $A \xcs B \xbe \xdf $

$(F3' )$ $A,B \xbe \xdf $ imply $A \xcs B \xEd \xCQ.$

So a weak filter satisfies $(F3' )$ instead of (F3).

An (weak) ideal on or over $X$ is a set $ \xdi \xcc \xdp (X),$ such that
$(I1)-(I3)$ $((I1),$ (I2),
$(I3' )$ respectively) hold:
(I1) $ \xCQ \xbe \xdi $

(I2) $A \xcc B \xcc X,$ $B \xbe \xdi $ imply $A \xbe \xdi $

(I3) $A,B \xbe \xdi $ imply $A \xcv B \xbe \xdi $

$(I3' )$ $A,B \xbe \xdi $ imply $A \xcv B \xEd X.$

So a weak ideal satisfies $(I3' )$ instead of (I3).

\ed

Elements of a filter on $X$ are called big subsets of $X,$ their
complements are
called small, and the rest have "medium size". The set of the
$X-$complements of
the elements of a filter form an ideal, and vice versa. Note that these
notions of "big" and "small" are by definition relative to the base set
$X.$ This
is a different, additional notion of relativity than the one we will see
in
Section 5.

These notions are related to nonmonotonic logics as follows:

We can say that, normally, $ \xbf $ implies $ \xbq $ iff in a big subset
of all $ \xbf -$cases, $ \xbq $
holds. In preferential terms, $ \xbf $ implies $ \xbq $ iff $ \xbq $ holds
in all minimal
$ \xbf -$models. If $ \xbm $ is the model choice function of a
preferential structure, i.e.
$ \xbm ( \xbf )$ is the set of minimal $ \xbf -$models, then $ \xbm ( \xbf
)$ will be a (the smallest)
big subset of the set of $ \xbf -$models, and the filter over the $ \xbf
-$models is the
pricipal filter generated by $ \xbm ( \xbf ).$

Due to the finite intersection property, filters and ideals work well with
logics:
If $ \xbf $ holds normally, as it holds in a big subset, and so does $
\xbf ',$ then
$ \xbf \xcu \xbf ' $ will normally hold, too, as the intersection of two
big subsets is
big again. This is a nice property, but not justified in all situations,
consider e.g. simple counting of a finite subset. (The question has a
name,
"lottery paradox": normally no single participant wins, but someone wins
in the
end.) This motivates the weak versions.

Normality defined by (weak or not) filters is a local concept: the filter
defined on $X$ and the one defined on $X' $ might be totally independent.
Consider,
however, the following two situations: Let $Y' $ be a big subset of $X',$
$X \xcc X',$
and $Y' \xcc X.$ If "size" has any absolute meaning, then $Y' $ should be
a big subset
of $X,$ too. On the other hand, let $X$ and $X' $ be big subsets of $Y,$
then
there are good reasons (analogue to those justifying the intersection
property
of filters) to assume that $X \xcs X' $ is also a big subset of $X'.$
These set
properties
are strongly connected to logical properties: For instance, if the latter
property holds, we can deduce the logical property Cautious Monotony (see
below for a formal definition): If $ \xbq $ implies normally $ \xbf $ and
$ \xbf ',$ because the
sets $X$ and $X' $ of $ \xbq \xcu \xbf -$models and $ \xbq \xcu \xbf '
-$models are big subsets of the set $Y$
of $ \xbq -$models, then $ \xbq \xcu \xbf ' $ will imply normally $ \xbf $
too, as the set $X \xcs X' $ of
$ \xbq \xcu \xbf \xcu \xbf ' -$models will be a big subset of the set $X'
$ of $ \xbq \xcu \xbf ' -$models.

More precisely, the reasoning is (a little simplified) as follows: If A
and $A' $
are small subsets of $B,$ then A will also be a small subset of B-A'.
Changing
the reference set $B$ just a little will not affect size. (In more detail:
If A and $A' $ are small subsets of $B,$ then A-A' will also be a small
subset of
$B-A'.)$ Note that "small" is used here in two conceptually very
different
ways: (1) $ \xbq \xcn \xbf $ iff the set of $ \xbq \xcu \xCN \xbf -$models
is a small subset of the
$ \xbq -$models, (2) when we change the reference class $B$ just a little
bit,
relative size is preserved.

Seen more
abstractly, such set properties allow the transfer of big subsets from one
to
another base set (and the conclusions drawn on this basis), and we call
them
"coherence properties". They are very important, not only for working with
a logic which respects them, but also for soundness and completeness
questions,
often they are at the core of such problems. The reader is invited to
read the articles by Ben-David and Ben-Eliyahu [BB94] and Friedman and
Halpern
[FH98], which treat essentially the same questions in different languages
(and
perhaps their comparison by the author in [Sch97-4] and [Sch04]).

These notions are tied to logical properties as follows:

The left hand column presents the single formula version, the center
column
the theory version (a theory is, for us, an arbitrary set of formulas),
the
right hand column the algebraic version, describing the choice function on
the
model set, e.g. $f(X) \xcc X$ corresponds to the rule $ \xbf \xcl \xbq $
implies $ \xbf \xcn \xbq $ in the
formula version, and to $ \ol{T} \xcc \ol{ \ol{T} }$ in the theory
version. A short discussion of some
of the properties follows the table.

(PR) is also called infinite conditionalization - we choose the name for
its
central role for preferential structures. Note that in the
presence of $( \xbm \xcc ),$ and if $ \xdy $ is closed under finite
intersections, $( \xbm PR)$
is equivalent to

$( \xbm PR' )$ $f(X) \xcs Y \xcc f(X \xcs Y).$

The system of rules (AND), (OR), (LLE), (RW), (SC), (CP), (CM), (CUM) is
also
called system $P$ (for preferential), adding (RM) gives the system $R$
(for
rationality or rankedness).

(LLE), (RW), (CCL) will all hold automatically, whenever we work with
fixed
model sets. (SC) corresponds to the choice of a subset. (CP) is somewhat
delicate, as it presupposes that the chosen model set is non-empty. This
might
fail in the presence of ever better choices, without ideal ones.
(PR) is an inifinitary version of one half of the deduction theorem: Let
$T$ stand
for $ \xbf,$ $T' $ for $ \xbq,$ and $ \xbf \xcu \xbq \xcn \xbs,$ so $
\xbf \xcn \xbq \xcp \xbs,$ but $( \xbq \xcp \xbs ) \xcu \xbq \xcl \xbs.$
(CUM) (whose most interesting half in our context is (CM)) may best be
seen as
normal use of lemmas: We have worked hard and found some lemmas. Now
we can take a rest, and come back again with our new lemmas. Adding them
to the
axioms will neither add new theorems, nor prevent old ones to hold.
(RM) is perhaps best understood by looking at big and small subsets. If
the
set of $ \xbf \xcu \xbq -$models is a big subset of the set of $ \xbf
-$models, and the
set of $ \xbf \xcu \xbq ' -$models is a not a small subset of the set of $
\xbf -$models (i.e.
big or of medium size), then the
set of $ \xbf \xcu \xbq \xcu \xbq ' -$models is a big subset of the set of
$ \xbf \xcu \xbq ' -$models.
Again, we have a double use of small/big here.

\bd

$\hspace{0.5em}$


\label{Definition 3.2}

{\small

\begin{tabular*}{11.5cm}{c@{\extracolsep\fill}|c|c}

(AND) \xEH (AND) \xEH \xEP

$ \xbf \xcn \xbq,  \xbf \xcn \xbq '   \xch $ \xEH
$ T \xcn \xbq, T \xcn \xbq '   \xch $ \xEH
\xEP

$ \xbf \xcn \xbq \xcu \xbq ' $ \xEH
$ T \xcn \xbq \xcu \xbq ' $ \xEH
\xEP

\hline

(OR) \xEH (OR) \xEH $( \xbm \xcv w)$ - $w$ for weak \xEP

$ \xbf \xcn \xbq,  \xbf ' \xcn \xbq   \xch $ \xEH
$T \xcn \xbq, T' \xcn \xbq   \xch $ \xEH
$f(A \xcv B) \xcc f(A) \xcv f(B)$
\xEP

$ \xbf \xco \xbf ' \xcn \xbq $ \xEH
$T \xco T' \xcn \xbq $ \xEH
\xEP

\hline

(LLE) or \xEH (LLE) \xEH \xEP

Left Logical Equivalence \xEH \xEH \xEP

$ \xcl \xbf \xcr \xbf ',  \xbf \xcn \xbq   \xch $ \xEH
$ \ol{T}= \ol{T' }  \xch   \ol{ \ol{T} }= \ol{ \ol{T' } }$ \xEH \xEP

$ \xbf ' \xcn \xbq $ \xEH \xEH \xEP

\hline

(RW) or Right Weakening \xEH (RW) \xEH \xEP

$ \xbf \xcn \xbq,  \xcl \xbq \xcp \xbq '   \xch $ \xEH
$ T \xcn \xbq,  \xcl \xbq \xcp \xbq '   \xch $ \xEH
\xEP

$ \xbf \xcn \xbq ' $ \xEH
$T \xcn \xbq ' $ \xEH
\xEP

\hline

(CCL) or Classical Closure \xEH (CCL) \xEH \xEP

\xEH $ \ol{ \ol{T} }$ is classically \xEH \xEP

\xEH closed \xEH \xEP

\hline

(SC) or Supraclassicality \xEH (SC) \xEH $( \xbm \xcc )$ \xEP

$ \xbf \xcl \xbq $ $ \xch $ $ \xbf \xcn \xbq $ \xEH $ \ol{T} \xcc \ol{
\ol{T} }$ \xEH $f(X) \xcc X$ \xEP
\hline

(CP) or \xEH (CP) \xEH $( \xbm \xCQ )$ \xEP

Consistency Preservation \xEH \xEH \xEP

$ \xbf \xcn \xcT $ $ \xch $ $ \xbf \xcl \xcT $ \xEH $T \xcn \xcT $ $ \xch
$ $T \xcl \xcT $ \xEH $f(X)= \xCQ $ $ \xch $ $X= \xCQ $ \xEP

\hline

(RM) or Rational Monotony \xEH (RM) \xEH $( \xbm =)$ \xEP

$ \xbf \xcn \xbq,  \xbf \xcN \xbq '   \xch $ \xEH
$T \xcn \xbq, T \xcN \xbq '   \xch $ \xEH
$X \xcc Y, Y \xcs f(X) \xEd \xCQ   \xch $
\xEP

$ \xbf \xcu \xbq ' \xcn \xbq $ \xEH
$T \xcv \{ \xbq ' \} \xcn \xbq $ \xEH
$f(X)=f(Y) \xcs X$ \xEP

\hline

(CM) or Cautious Monotony \xEH (CM) \xEH \xEP

$ \xbf \xcn \xbq,  \xbf \xcn \xbq '   \xch $ \xEH
$T \xcc \ol{T' } \xcc \ol{ \ol{T} }  \xch $ \xEH
$f(X) \xcc Y \xcc X  \xch $
\xEP

$ \xbf \xcu \xbq \xcn \xbq ' $ \xEH
$ \ol{ \ol{T} } \xcc \ol{ \ol{T' } }$ \xEH
$f(Y) \xcc f(X)$
\xEP

\hline

(CUM) or Cumulativity \xEH (CUM) \xEH $( \xbm CUM)$ \xEP

$ \xbf \xcn \xbq   \xch $ \xEH
$T \xcc \ol{T' } \xcc \ol{ \ol{T} }  \xch $ \xEH
$f(X) \xcc Y \xcc X  \xch $
\xEP

$( \xbf \xcn \xbq '   \xcj   \xbf \xcu \xbq \xcn \xbq ' )$ \xEH
$ \ol{ \ol{T} }= \ol{ \ol{T' } }$ \xEH
$f(Y)=f(X)$ \xEP

\hline

\xEH (PR) \xEH $( \xbm PR)$ \xEP

$ \ol{ \ol{ \xbf \xcu \xbf ' } }$ $ \xcc $ $ \ol{ \ol{ \ol{ \xbf } } \xcv
\{ \xbf ' \}}$ \xEH
$ \ol{ \ol{T \xcv T' } }$ $ \xcc $ $ \ol{ \ol{ \ol{T} } \xcv T' }$ \xEH
$X \xcc Y$ $ \xch $
\xEP

\xEH \xEH $f(Y) \xcs X \xcc f(X)$
\xEP

\end{tabular*}

}

\ed

We recall that smooth preferential models generate the system $P,$ whereas
smooth ranked models generate the system $R$ - see e.g.  \cite{KLM90} or
\cite{Sch04}
for details.
\section{
Interpretations
}
\label{Section 4}
\subsection{
Introduction
}
\label{Section 4.1}

We will discuss in this Section three interpretations of inheritance nets.

First, we will indicate fundamental differences between inheritance and
the
systems $P$ and $R,$ they will be elaborated in Section 5, where an
interpretation
in terms of small sets will be tried nonetheless.

Second, we will interpret inheritance nets as systems of information and
information flow.

Third, we will interpret inheritance nets as systems of prototypes.

Inheritance nets present many intuitively attractive properties, thus it
is
not surprising that we can interpret them in several ways. Similarly,
preferential structures can be used as a semantics for deontic, and for
nonmonotonic logic, they express a common idea: choosing a subset of
models
by a binary relation. Thus, such an ambiguity need not be a sign for a
basic
flaw.
\subsection{
Informal comparison of inheritance with the systems P and R
}
\label{Section 4.2}

\paragraph{
The main issues
}

$\hspace{0.01em}$

In the author's opinion, the following two properties of inheritance
diagrams
show the deepest difference to preferential and similar semantics, and the
first even to classical logic. They have to be taken seriously, as they
are at
the core of inheritance systems, are independent of the particular
formalism,
and show that there is a fundamental difference between the former and the
latter. Consequently, any attempt at translation will have to stretch one
or
both sides perhaps beyond the breaking point.

(1) Relevance,

(2) subideal situations, or relative normality

Both (and more) can be illustrated by the following simple Diagram \ref{Diagram
4.1}
(which also shows conflict resolution by specificity).

(1) Relevance: As there is no monotonous path whatever between $e$ and
$d,$
the question whether e's are d's or not, or vice versa, does not even
arise.
For the same reason, there is no question whether b's are $c' $s, or not.
(As
a matter of fact, we will see below that b's are non-c's in system $P$ -
see
Definition \ref{Definition 3.2}).
In upward chaining formalisms, as there is no valid positive path from a
to $d,$
there is no question either whether a's are f's or not.

The problem of relevance has a trivial answer in upward chaining
inheritance
nets: there must be some valid initial segment in the diagram to establish
relevance.

Of course, in classical logic, all information is relevant to the rest, so
we
can say e.g. that e's are $d' $s, or e's are $non-d' $s, or some are $d'
$s, some are
not, but there is a connection. As preferential models are based on
classical
logic, the same argument applies to them.

(2) In our diagram, a's are $b' $s, but not ideal $b' $s, as they are not
$d' $s, the
more specific information from $c$ wins. But they are $e' $s, as ideal b's
are.
So they are not perfectly ideal $b' $s, but as ideal b's as possible.
Thus, we
have graded ideality, which does not exist in preferential and similar
structures. In those structures, if an element is an ideal element, it has
all
properties of such,
if one such property is lacking, it is not ideal, and we $can' t$ say
anything any
more. Here, however, we sacrifice as little normality as possible, it is
thus a minimal change formalism.

In comparison, questions of information transfer and strength of
information
seem lesser differences. Already systems $P$ and $R$ (see Definition
\ref{Definition 3.2})
differ on information transfer. In both cases, transfer is
based on the same notion of smallness, which describes ideal situations.
But,
as said in Section \ref{Section 3}, conceptually, this is very different
from the use of
"ideal-smallness", describing normal situations. Thus,
it can be considered also on this level an independent question, and we
can
imagine systems based on absolutely ideal situations for normality, but
with a
totally different transfer mechanism.

\bfc

\bdia

\label{Diagram 4.1}

\unitlength1.0mm
\begin{picture}(130,110)

\newsavebox{\Tweety}
\savebox{\Tweety}(130,110)[bl]
{

\put(57,18){\vector(1,1){24}}
\put(51,18){\vector(-1,1){24}}
\put(27,51){\vector(1,1){24}}
\put(81,51){\vector(-1,1){24}}

\put(67,61){\line(1,1){4}}

\put(81,47){\vector(-1,0){54}}

\put(24,51){\vector(0,1){22}}
\put(54,81){\vector(0,1){22}}

\put(53,16){$a$}
\put(23,46){$b$}
\put(83,46){$c$}
\put(53,76){$d$}
\put(23,76){$e$}
\put(53,106){$f$}

}

\put(0,0){\usebox{\Tweety}}
\end{picture}

\edia

\efc

For these reasons, extending preferential and related semantics to cover
inheritance nets seems to stretch them to the breaking point, Thus, we
should
also look for other interpretations.
(The term "interpretation" is used here in a non-technical sense.) In
particular, it seems worth while to connect inheritance systems to other
problems, and see whether there are similarities there. This is what we do
now.
We come back to the systems $P$ and $R$ in Section 5.
\subsection{
Inheritance as information transfer
}
\label{Section 4.3}

An informal argument showing parallel ideas common to inheritance with an
upward chaining formalism and information transfer is as follows: First,
arrows represent certainly some kind of information, of the kind "most a's
are
$b' $s" or so. (See Diagram \ref{Diagram 4.1}.) Second, to be able to
use information,
e.g. "d's are $f' $s" at a, we have
to be able to connect from a to $d$ by a valid path, this information has
to be
made accessible to a, or, in other terms, a working information channel
from
a to $d$ has to be established. Third, specificity (when present) decides
conflicts (we take the split validity approach). This can be done
procedurally, or, perhaps simpler and certainly in a more transparent way,
by assigning a comparison of information strength to valid paths. Now,
information strength may also be called truth value (to use a term
familiar
in logic) and the natural entity at hand is the node itself - this is just
a
cheap formal trick without any conceptual meaning.

When we adopt this view, nodes and
arrows (and valid paths) have multiple functions, and it may seem that we
overload the (deceptively) simple picture. But it is perhaps the charm and
the utility and naturalness of inheritance systems that they are not
"clean",
and hide many complications under a simple surface, as human common sense
reasoning often does, too.

In a certain way, this is a poor man's interpretation, as it does not base
inheritance on another formalism, but gives only an intuitive reading.
Yet, it
gives a connection to other branches of reasoning, and is as such already
justified - in the author's opinion. Moreover, our analysis makes a clear
distinction between arrows and composite valid paths. This distinction
is implicit in inheritance formalisms, we make it explicit through our
concepts.

But this interpretation is by no means the
only one, and can only be suggested as a possibility.

We will now first give the details, and then discuss our interpretation.

\paragraph{
(1) Information:
}

$\hspace{0.01em}$

Direct arrows (negative or positive) represent information, valid for
their
source. Thus, in a set reading, if there is an arrow $A \xcp B$ in the
diagram, most
elements of $ \xCB $ will be in $B,$ in short: "most $ \xCB ' $s are $B'
$s" - and $A \xcP B$ will mean
that most $ \xCB ' $s are not $B' $s.

\paragraph{
(2) Information sources and flow:
}

$\hspace{0.01em}$

Nodes are information sources. If $A \xcp B$ is in the diagram, $ \xCB $
is the source
of the information "most $ \xCB ' $s are $B' $s".

A valid, composed or atomic positive path $ \xbs $ from $U$ to $ \xCB $
makes the
information of source $ \xCB $ accessible to $U.$ One might also say that
$ \xCB ' $s
information becomes relevant to $U.$ Otherwise, information is considered
independent - only (valid) positive paths create the dependencies.

(If we want to conform to inheritance, we must not add trivialities like
"x's are $x' $s", as this would require $x \xcp x$ in the corresponding
net, which, of
course, will not be there in an acyclic net.)

\paragraph{
(3) Information strength:
}

$\hspace{0.01em}$

A valid, composed or atomic positive path $ \xbs $ from $ \xCB ' $ to $
\xCB $ allows us to compare
the strength of information source $ \xCB ' $ with that of $ \xCB:$ $
\xCB ' $ is stronger than $ \xCB.$
(In the set reading, this comparison is the result of specificity: more
specific
information is considered more reliable.) If there is no such valid path,
we
cannot resolve contradictions between information from $ \xCB $ and $ \xCB
'.$
This interpretation results in split validity preclusion: the
comparison between information sources $ \xCB ' $ and $ \xCB $ is
absolute, and does NOT
depend on the $U$ from which both may be accessible - as can be the case
with
total validity preclusion. Of course, if desired, we can also adopt this
(much more complicated) idea.

Nodes are also truth values. They are the strength of the information
whose
source they are. (This might seem an abuse of nodes, but we already have
them,
so why not use them?)

\paragraph{
Discussion:
}

$\hspace{0.01em}$

Considering direct arrows as information meets probably with little
objection.

The conclusion of a valid path (e.g. if $ \xbs:a \xcp b$ is valid, then
its conclusion
is "a's are $b' $s") is certainly also information, but it has a status
different
from the information of a direct link, so we should distinguish it
clearly. At
least in upward chaining formalisms, using the path itself as some channel
through which information flows, and not the conclusion, seems more
natural.
The conclusion says little about the inner structure of the path, which is
very important in inheritance networks, e.g. for preclusion. When
calculating
validity
of paths, we look at (sub- and other) paths, but not their results, and
should
also express this clearly.

Once we accept this picture of valid positive paths as information
channels, it
is natural to see their upper ends as information sources.

Our interpretation supports upward chaining, and vice versa, upward
chaining
supports our interpretation.

One of the central ideas of inheritance is preclusion, which, in the case
of
split validity preclusion, works by an absolute comparison between nodes.
Thus,
if we accept split validity preclusion, it is natural to see valid
positive
paths as
comparisons between information of different strengths. Conversely, if we
accept absolute comparison of information, we should also accept split
validity
preclusion - these interpretations support each other.

Whatever type of preclusion we accept, preclusion clearly compares
information
strength, and allows us to decide for the stronger one. We can see this
procedurally, or by giving different values to different information
(depending
on the source), which we can call truth values to connect our picture to
other
areas of logic. It is then natural - as we have it already - to use the
source
node itself as truth value, with comparison via valid positive paths.

\paragraph{
Illustration:
}

$\hspace{0.01em}$

Thus, in a given node $U,$ information from $ \xCB $ is accessible iff
there is a
valid positive path from $U$ to $ \xCB,$ and if information from $ \xCB '
$ is also accessible,
and there is a valid positive path from $ \xCB ' $ to $ \xCB,$ then, in
case of conflict,
information from $ \xCB ' $ wins over that from $ \xCB,$ as $ \xCB ' $
has a better truth value.
In the Tweety diagram, see Diagram \ref{Diagram 2.1},
Tweety has access to penguins and birds, the horizontal link from penguin
to
bird compares the strengths, and the fly/not fly arrows are the
information.

A negative direct link can only be information. A positive direct
link is information at its source, but it can also be a comparison of
truth
values, or it can give access from its source to information at its end.
A valid positive, composed path can only be comparison of truth values, or
give
access to information, it is NOT information itself. A valid negative
composed
path has no function, only its parts have.

We obtain automatically that direct information is
stronger than any other information: If $ \xCB $ has information $ \xbf,$
and there is a
valid path from $ \xCB $ to $B,$ making B's information accessible to $
\xCB,$ then this same
path also compares strength, and $ \xCB ' $s information is stronger than
B's
information. Seen from $ \xCB,$ i.e. just considering information
accessible to $ \xCB,$ $ \xCB ' $s
own information will always be best.

Our interpretation underlines the importance of initial segments: Initial
segments make information accessible. Thus, initial segments have to be
valid.

Inheritance diagrams in this interpretation do not only
represent reasoning with many truth values, but also reasoning ABOUT those
truth
values: their comparison is done by the same underlying mechanism.

\paragraph{
Further comments:
}

$\hspace{0.01em}$

Our reading also covers enriched diagrams, where arbitrary
information can be "appended" to a node.

An alternative way to see a source of information is to see it as a reason
to believe the information it gives. $U$ needs a reason to believe
something, i.e. a
valid path from $U$ to the source of the information, and also a reason to
disbelieve, i.e. if $U' $ is below $U,$ and $U$ believes and $U' $ does
NOT believe some
information of $ \xCB,$ then either $U' $ has
stronger information to the contrary, or there is not a valid path to $
\xCB $ any more
(and neither to any other possible source of this information).
$('' Reason'',$ a concept very important in this context, was introduced
by A.Bochman
into the discussion.)

The restriction that negative links can only be information applies to
traditional inheritance networks, and the
author makes no claim whatever that it should also hold for modified such
systems, or in still other contexts. One of the reasons why we do not have
"negative nodes", and thus negated arrows also in the middle of paths
might be
the following (with $ \xdC $ complementation): If, for some $X,$ we also
have a node for
$ \xdC X,$ then we should have $X \xcP \xdC X$ and $ \xdC X \xcP X,$ thus
a cycle, and arrows from $Y$ to
$X$ should be accompanied by their opposite to $ \xdC X,$ etc. This would
complicate
the picture, perhaps without any real gain in insight.

We translate the analysis and decision of Definition \ref{Definition 2.2}
now into the
picture of information sources, accessibility, and comparison via valid
paths.
This is straightforward:

(1) We have that information from $A_{i},$ $i \xbe I,$ about $B$ is
accessible from $U,$ i.e.
there are valid positive paths from $U$ to all $A_{i}.$ Some $A_{i}$ may
say $ \xCN B,$ some $B.$

(2) If information from $A_{i}$ is comparable with information from
$A_{j}$ (i.e. there
is a valid positive path from $A_{i}$ to $A_{j}$ or the other way around),
and $A_{i}$
contradicts $A_{j}$ with respect to $B,$ then the weaker information is
discarded.

(3) There remains a (nonempty, by lack of cycles) set of the $A_{i},$ such
that
for no such
$A_{i}$ there is $A_{j}$ with better contradictory information about $B.$
If the information
from this remaining set is contradictory, we accept none (and none of the
paths
either), if not, we accept the common conclusion and all these paths.

We continue now Remark \ref{Remark 2.1}, (4), and turn this into a
formal system.

Fix a diagram $ \xbG,$ and do an induction as in Definition \ref{Definition
2.2}.

\bd

$\hspace{0.01em}$

(+++*** Orig. No.:  Definition 4.1 )

\label{Definition 4.1}

(1) We distinguish $a \xch b$ and $a \xch_{x}b,$ where the intuition of $a
\xch_{x}b$ is:
we know with strength $x$ that a's are $b' $s, and of $a \xch b$ that it
has been
decided taking all information into consideration that $a \xch b$ holds.

(2) $a \xcp b$ implies $a \xcp_{a}b,$ likewise $a \xcP b$ implies $a
\xcp_{a} \xCN b.$

(3) $a \xch_{a}b$ implies $a \xch b,$ likewise $a \xch_{a} \xCN b$ implies
$a \xch \xCN b.$ This expresses
the fact that direct arrows are uncontested.

(4) $a \xch b$ and $b \xch_{b}c$ imply $a \xch_{b}c,$ likewise for $b
\xch_{b} \xCN c.$ This expresses
concatenation - but without deciding if it is accepted! Note we cannot
make
$a \xch b$ and $b \xch c$ imply $a \xch_{b}c$ a rule, as this would make
concatenation of two
composed paths possible.

(5) We decide acceptance of composed paths as in Definition \ref{Definition
2.2}, where
preclusion uses accepted paths for deciding.

\ed

Note that we also reason in this system about relative strength of truth
values,
which are just nodes, this is then, of course, used in the acceptance
condition,
in preclusion, more precisely.
\subsection{
Inheritance as reasoning with prototypes
}
\label{Section 4.4}

Some of the issues we discuss here apply also to the more general picture
of
information and its transfer. We present them here for motivational
reasons:
it seems easier to discuss them in the (somewhat!) more concrete setting
of prototypes than in the very general situation of information handling.
These issues will be indicated.

It seems natural to see information in inheritance networks as information
about prototypes. (We do not claim that our use of the word "prototype"
has
more than a vague relation to the use in psychology. We do not try to
explain
the usefulness of prototypes either, one possibility is that there are
reasons
why birds fly, and why penguins $don' t,$ etc.) In the Tweety diagram, we
will thus say that prototypical birds will fly, prototypical penguins will
not
fly. More precisely, the property "fly" is part of the bird prototype, the
property $'' \xCN fly'' $ part of the penguin prototype. Thus, the
information is
given for some node, which defines its application or domain (bird or
penguin in
our example) - beyond this node, the property is not defined (unless
inherited,
of course). It might very well be that no element of the domain has ALL
the
properties of the prototype, every bird may be exceptional in some sense.
This
again shows that we are very far from the ideal picture of small and big
subsets as used in systems $P$ and $R.$ (This, of course, goes beyond the
problem
of prototypes.)

Of course, we will want to "inherit" properties of prototypes,
for instance in Diagram \ref{Diagram 4.1},
a "should" inherit the property $e$ from $b,$ and
the property $ \xCN d$ from $c.$ Informally, we will argue as follows:
"Prototypical
a's have property $b,$ and prototypical b's have property $e,$ so it seems
reasonable to assume that prototypical a's also have property $e$ - unless
there
is better information to the contrary." A plausible solution is then to
use
upward chaining inheritance as described above to find all relevant
information, and then compose the prototype.

We discuss now three points whose importance goes beyond the treatment of
prototypes:

(1) Using upward chaining has an
additional intuitive appeal: We consider information
at a the best, so we begin with $b$ (and $c),$ and only then, tentatively,
add
information $e$ from $b.$ Thus, we begin with strongest information, and
add weaker
information successively - this seems good reasoning policy.

(2) In upward chaining, we also collect information
at the source (the end of the path), and do not use information which was
already filtered by going down - thus the information we collect has no
history, and we cannot encounter problems of iterated revision, which are
problems of history of change. (In downward chaining, we
only store the reasons why something holds, but not why something does not
hold,
so we cannot erase this negative information when the reason is not valid
any
more. This is an asymmetry apparently not much noted before.
Consider Diagram \ref{Diagram 2.4}. Here, the
reason why $u$ does not accept $y$ as information, but $ \xCN y,$ is the
preclusion
via $x.$ But from $z,$ this preclusion is not valid any more, so the
reason
why $y$ was rejected is not valid any more, and $y$ can now be accepted.)

(3) We come back to the question of extensions vs. direct scepticism.
Consider the Nixon Diamond, Diagram \ref{Diagram 2.2}. Suppose Nixon
were a subclass
of Republican and Quaker. Then the extensions approach reasons as follows:
Either the Nixon class prototype has the pacifist property, or the hawk
property, and we consider these two possibilities. But this is not
sufficient:
The Nixon class prototype might have neither property - they are normally
neither pacifists, nor hawks, but some are this, some are that. So the
conceptual basis for the extensions approach does not hold: "Tertium non
datur"
just simply does not hold (as in Intuitionist Logic, where we may have
neither
a proof for $ \xbf,$ nor for $ \xCN \xbf ).$

Once we fixed this decision, i.e. how to find the relevant information, we
can
still look upward or downward in the net and investigate the changes
between the
prototypes in going upward or downward,
as follows: E.g., in above example, we can look at the node a and its
prototype,
and then at the change going from a to $b,$ or, conversely, look at $b$
and its
prototype, and then at the
change going from $b$ to a. The problem of finding the information, and
this
dynamics of information change have to be clearly separated.

In both cases, we see the following:

(1) The language is kept small, and thus efficient.

For instance, when we go from a to $b,$ information about $c$ is lost, and
"c" does
not figure any more in the language, but $f$ is added. When we go from $b$
to a,
$f$ is lost, and $c$ is won. In our simple picture, information is
independent,
and contradictions are always between two bits of information.

(2) Changes are kept small, and need a reason to be effective.
Contradictory, stronger information will override the old one, but
no other information, except in the following case:
making new information (in-) accessible will cause indirect changes,
i.e. information now made (in-) accessible via the new node.
This is similar to formalisms of causation: if a reason is not there
any more, its effects vanish, too.

It is perhaps more natural when going downward also to
consider "subsets", as follows: Consider Diagram \ref{Diagram 4.1}.
b's are $d' $s, and c's are $ \xCN d' $s, and c's are also $b' $s. So it
seems
plausible to go beyond the language of inheritance nets, and conclude that
b's which are not c's will be $d' $s, in short to consider $(b-c)' $s. It
is obvious
which such subsets to consider, and how to handle them: For instance,
loosely
speaking, in $b \xcs d$ $e$ will hold, in $b \xcs c \xcs d$ $ \xCN f$ will
hold, in $b \xcs d \xcs \xdC c$ $f$ will
hold, etc. This is just putting the bits of information together.

We turn to another consideration, which will also transcend the prototype
situation and we will (partly) use the intuition that nodes stand for
sets, and
arrows for (soft) inclusion in a set or its complement.

In this reading, specificity stands for soft set inclusion. So, if $b$ and
$c$ are
visible from a, and there is a valid path from $c$ to $b$ (as in Diagram
\ref{Diagram 4.1}),
then a is a subset both of $b$ and $c,$ and $c$ a subset of $b,$ so
$a \xcc c \xcc b$ (softly). But then a is closer to $c$ than a is to $b.$
Automatically,
a will be closest to itself.

When we go now from $b$ to $c,$ we loose information $d$ and $f,$ win
information
$ \xCN d,$ but keep information $e.$ Thus, this is minimal change: we give
up (and
win) only the necessary information, but keep the rest. As our language is
very simple, we can use the Hamming distance between formula sets here.
(We
will make a remark on more general situations just below.)

When we look now again from a, we take the set-closest class (c), and use
the information of $c,$ which was won by minimal change (i.e. the Hamming
closest)
from information of $b.$ So we have the interplay of two distances, where
the set distance certainly is not symmetrical, as we need valid paths for
access and comparison. If there is no such valid path, it is reasonable to
make the distance infinite.

The promised remark on more general situations: in richer languages, we
cannot
count formulas to determine the Hamming distance between two situations
(i.e.
models or model sets), but have to take the difference in propositional
variables. Consider e.g. the language with two variables, $p$ and $q.$ The
models
(described by) $p \xcu q$ and $p \xcu \xCN q$ have distance 1, whereas $p
\xcu q$ and $ \xCN p \xcu \xCN q$
have distance 2. Note that this distance is NOT robust under re-definition
of the language. Let $p' $ stand for $(p \xcu q) \xco ( \xCN p \xcu \xCN
q),$ and $q' $ for $q.$ Of course,
$p' $ and $q' $ are equivalent descriptions of the same model set, as we
can define
all the old singletons also in the new language. Then
the situations $p \xcu q$ and $ \xCN p \xcu \xCN q$ have now distance 1,
as one corresponds to
$p' \xcu q',$ the other to $p' \xcu \xCN q'.$

There might be misunderstandings about the use of the word "distance"
here. The
author is fully aware that inheritance networks cannot be captured by
distance
semantics in the sense of preferential structures. But we do NOT think
here of
distances from one fixed ideal point, but of relativized distances: Every
prototype is the origin of measurements. E.g., the bird prototype is
defined
by "flying, laying eggs, having feathers  \Xl.". So we presume that all
birds
have these properties of the prototype, i.e. distance 0 from the
prototype.
When we see that penguins do not fly, we move as little as possible from
the
bird prototype, so we give up "flying", but not the rest. Thus, penguins
(better: the penguin prototype) will have distance 1 from the bird
prototype
(just one property has changed). So there is a new prototype for penguins,
and
considering penguins, we will not measure from the bird prototype, but
from the
penguin prototype, so the point of reference changes. This is exactly as
in
distance semantics for theory revision, introduced in  \cite{LMS01},
only the point
of reference is not the old theory $T,$ but the old prototype, and the
distance
is a very special one, counting properties assumed to be independent. (The
picture is a little bit more complicated, as the loss of one property
(flying)
may cause other modifications, but the simple picture suffices for this
informal argument.)

We conclude this Section with two remarks, the first on prototypes, the
second a general remark on preclusion in inheritance systems.

Realistic prototypical reasoning will probably neither always be upward,
nor
always be downward. A medical doctor will not begin with the patient's
most
specific category (name and birthday or so), nor will he begin with all he
knows about general objects. Therefore, it seems reasonable to investigate
upward and downward reasoning here.

Obviously, in some cases, it need not be specificity, which decides
conflicts. Consider the case where Tweety is a bird, but a dead animal.
Obviously, Tweety will not fly, here because the predicate "dead" is very
strong and overrules many normal properties. When we generalize this, we
might
have a hierarchy of causes, where one overrules the other, or the result
may be
undecided. For instance, a falling object might be attracted in a magnetic
field, but a gusty wind might prevent this, sometimes, with unpredictable
results. This is then additional information (strength of cause), and this
problem is not addressed directly in traditional inheritance networks, we
would have to introduce a subclass "dead bird" - and subclasses often have
properties of "pseudo-causes", being a penguin probably is not a "cause"
for not flying, nor bird for flying, still, things change from class to
subclass for a reason.
\section{
Detailed translation of inheritance to modified systems of small sets
}
\label{Section 5}
\subsection{
Normality
}
\label{Section 5.1}

As we saw already in Section \ref{Section 4.2},
normality in inheritance (and Reiter defaults etc.) is relative, and as
much normality as possible is preserved. There is no absolute $N(X),$ but
only
$N(X, \xbf ),$ and $N(X, \xbf )$ might be defined, but not $N(X, \xbq ).$
Normality in the sense of preferential structures is absolute: if $x$ is
not
in $N(X)$ $(=$ $ \xbm (X)$ in preferential reading), we do not know
anything beyond
classical logic.
This is the dark Swedes' problem: even dark Swedes should probably be
tall.
Inheritance systems are different: If birds usually lay eggs, then
penguins,
though abnormal with respect to flying, will still usually lay eggs.
Penguins are fly-abnormal birds, but will continue to be egg-normal birds
-
unless we have again information to the contrary.
So the absolute, simple $N(X)$ of preferential structures splits up into
many, by default independent, normalities, $N(X, \xbf )$ for $ \xbf
-$normal etc.
This corresponds to intuition: There are no absolutely normal birds, each
one is particular in some sense, so $ \xcS \{N(X, \xbf ): \xbf \xbe \xdl
\}$ may well be empty, even
if each single $N(X, \xbf )$ is almost all birds.

What are the laws of relative normality?
$N(X, \xbf )$ and $N(X, \xbq )$ will be largely independent (except for
trivial situations,
where $ \xbf \xcr \xbq,$ $ \xbf $ is a tautology, etc.). $N(X, \xbf )$
might be defined, and $N(X, \xbq )$
not. Thus, if there is no arrow, or no path, between $X$ and $Y,$ then
$N(X,Y)$ and
$N(Y,X)$ - where X,Y are also properties - need not be defined. This will
get rid
of the unwanted connections found with absolute normalities, as
illustrated
by Fact \ref{Fact 5.1}.

We interpret now "normal" by "big set", i.e. essentially $'' \xbf $ holds
normally in
$X$ iff there is a big subset of $X,$ where $ \xbf $ holds". This will, of
course, be
modified.
\subsection{
Small sets
}
\label{Section 5.2}

The main interest of this Section is perhaps to show the adaptations of
the
concept of small/big subsets necessary for a more "real life" situation,
where we have to relativize. The amount of changes we do illustrates the
problems and what can be done, but also perhaps what should not be done,
as the
concept is stretched too far.

The usual informal way of speaking about inheritance networks
(plus other considerations) motivates an interpretation by
sets and soft set inclusion - $A \xcp B$ means that "most $ \xCB ' $s are
$B' $s". Just as with
normality, the "most" will have to be relativized, i.e. there is a
$B-$normal part
of $ \xCB,$ and a $B-$abnormal one, and the first is $B-$bigger than the
second - where
"bigger" is relative to $B,$ too. (A further motivation for this set
interpretation
is the often evoked specificity argument for preclusion. Thus, we will now
translate our remarks about normality into the language of big and small
subsets.)

Recall our remarks about relative normality. $N(X, \xbf )$ is, a priori,
independent of $N(X, \xbq ),$ and $N(X, \xbf )$ might be defined, but not
$N(X, \xbq ).$ Thus, we
will have $ \xbf -$big subsets of $X,$ and $ \xbq -$big subsets, and the
two are essentially
independent, may perhaps even have an empty intersection, only one may be
defined, etc.

Consider now the system $P$ (with Cumulativity), see Definition \ref{Definition
3.2}.
Small sets (see Definition \ref{Definition 3.1}) are used in two
conceptually very distinct ways: $ \xba \xcn \xbb $ iff the set of $ \xba
\xcu \xCN \xbb -$cases is a small
subset (in the absolute sense, there is just one system of big subsets of
the
$ \xba -$cases) of the set of $ \xba -$cases. The second use is in
information transfer,
used in Cumulativity, or Cautious Monotony more precisely: if the set of
$ \xba \xcu \xCN \xbg -$cases is a small
subset of the set of $ \xba -$cases, then $ \xba \xcn \xbb $ carries over
to $ \xba \xcu \xbg:$ $ \xba \xcu \xbg \xcn \xbb.$
(See also the discussion in  \cite{Sch04}, page 86, after Definition
2.3.6.) It is
this transfer which we will consider here, and not things like AND, which
connect different $N(X, \xbf )$ for different $ \xbf.$

Before we go into details, we will show that e.g. the system $P$ is too
strong
to model inheritance systems, and that e.g. the system $R$ is to weak for
this
purpose. Thus, preferential systems are really quite different from
inheritance
systems.

\bfa

$\hspace{0.01em}$

(+++*** Orig. No.:  Fact 5.1 )

\label{Fact 5.1}

(a) System $P$ is too strong to capture inheritance.

(b) System $R$ is too weak to capture inheritance.

\efa

\paragraph{
Proof:
}

$\hspace{0.01em}$

(a) Consider the Tweety diagram, Diagram \ref{Diagram 2.1}. $c \xcp b
\xcp d,$ $c \xcP d.$ There is no
arrow $b \xcP c,$ and we will see that $P$ forces one to be there. For
this, we take
the natural translation, i.e. $X \xcp Y$ will be $'' X \xcs Y$ is a big
subset of $X'',$ etc. We
show that $c \xcs b$ is a small subset of $b,$ which we write $c \xcs
b<b.$
$c \xcs b=(c \xcs b \xcs d) \xcv (c \xcs b \xcs \xdC d).$ $c \xcs b \xcs
\xdC d \xcc b \xcs \xdC d<b,$ the latter by $b \xcp d,$ thus
$c \xcs b \xcs \xdC d<b,$ essentially by Right Weakening. Set now $X:=c
\xcs b \xcs d.$ As $c \xcP d,$
$X:=c \xcs b \xcs d \xcc c \xcs d<c,$ and by the same reasoning as above
$X<c.$ It remains to show
$X<b.$ We use now $c \xcp b.$ As $c \xcs \xdC b<c,$ and $c \xcs X<c,$ by
Cumulativity $X=c \xcs X \xcs b<c \xcs b,$
so essentially by OR $X=c \xcs X \xcs b<b.$ Using the filter property, we
see that
$c \xcs b<b.$

(b) Second, even $R$ is too weak: In the diagram $X \xcp Y \xcp Z,$ we
want to conclude that
most of $X$ is in $Z,$ but as $X$ might also be a small subset of $Y,$ we
cannot
transfer the information "most Y's are in $Z'' $ to $X.$

$ \xcz $
\\[3ex]

We have to distinguish direct information/arrows from inherited
information/valid paths. In the language of big/small sets, it is easiest
to do this by two types of big subsets: big ones and very big ones. We
will
denote the first big, the second BIG. This corresponds to the distinction
between $a \xch b$ and $a \xch_{a}b$ in Definition \ref{Definition 4.1}.

In particular, we will have the implications $BIG \xcp big$ and $SMALL
\xcp small,$ so
we have nested systems. Such systems were discussed in  \cite{Sch95-1}, see also
 \cite{Sch97-2}.

In particular, this distinction seems to be necessary to prevent arbitrary
concatenation of valid paths to valid paths, which would lead to
contradictions.
Consider e.g. $a \xcp b \xcp c \xcp d,$ $a \xcp e \xcP d,$ $e \xcp c.$
Then concatenating $a \xcp b$ with $b \xcp c \xcp d,$
both valid, would lead to a simple contradiction with $a \xcp e \xcP d,$
and not to
preclusion, as it should be - see below.

For the situation $X \xcp Y \xcp Z,$ we will then conclude that:

If $Y \xcs Z$ is a $Z-$BIG subset of $Y$ and $X \xcs Y$ is a $Y-$big
subset of $X$ then $X \xcs Z$ is a
$Z-$big subset of $X.$ (We generalize already to the case where there is a
valid
path from $X$ to $Y.)$

We call this procedure information transfer.

$Y \xcp Z$ expresses the direct information in this context, so $Y \xcs Z$
has to be a $Z-$BIG
subset of $Y.$

$X \xcp Y$ can be direct information, but it is used here as channel of
information
flow, in particular it might be a composite valid path, so in our context,
$X \xcs Y$ is a $Y-$big subset of $X.$

$X \xcs Z$ is a $Z-$big subset of $X:$ this can only be big, and not BIG,
as we have a
composite path.

The translation into big/small subsets and their modifications is now
quite
complicated: we seem to have to relativize, and we seem to need two types
of big/small. This casts, of course, a doubt on the enterprise of
translation.
The future will tell if any of the ideas can be used in other contexts.

We investigate this situation now in more detail, first without conflicts.

The way we cut the problem is not the only possible one. We were guided by
the
idea that we should stay close to usual argumentation about big/small
sets,
should proceed carefully, i.e. step by step, and should take a general
approach.

Note that we start without any $X-$big subsets defined, so $X$ is not even
a $X-$big
subset of itself.

(A) The simple case of two arrows, and no conflicts.

(In slight generalization:) If information $ \xbf $ is appended at $Y,$
and $Y$ is
accessible from $X$ (and there is no better information about $ \xbf $
available), $ \xbf $
will be valid at $X.$ For simplicity, suppose there is a direct positive
link from
$X$ to $Y,$ written sloppily $X \xcp Y \xcm \xbf.$

In the big subset reading, we will interpret this as: $Y \xcu \xbf $ is a
$ \xbf -$BIG
subset of $Y.$ It is important that this is now direct information, so we
have
"BIG" and not "big".

We read now $X \xcp Y$ also as: $X \xcs Y$ is an $Y-$big subset of $X$ -
this is the channel,
so just "big".

We want to conclude by transfer that $X \xcs \xbf $ is a $ \xbf -$big
subset of $X.$

We do this in two
steps: First, we conclude that $X \xcs Y \xcs \xbf $ is a $ \xbf -$big
subset of $X \xcs Y,$ and then, as
$X \xcs Y$ is an $Y-$big subset of $X,$ $X \xcs \xbf $ itself is a $ \xbf
-$big subset of $X.$ We do NOT
conclude that $(X-Y) \xcs \xbf $ is a $ \xbf -$big subset of X-Y, this is
very important, as we
want to preserve the reason of being $ \xbf -$big subsets - and this goes
via $Y!$

The transition from "BIG" to "big" should be at the first step, where we
conclude that $X \xcs Y \xcs \xbf $ is a $ \xbf -$big (and not $ \xbf
-$BIG) subset of $X \xcs Y,$ as it is
really here where things happen, i.e. transfer of information from $Y$ to
arbitrary subsets $X \xcs Y.$

Now, for the two steps in a slightly modified notation, corresponding to
the
diagram $X \xcp Y \xcp Z:$

(Here and in what follows, we will be very cautious, and relativize all
normalities. We could perhaps obtain our objective with a more daring
approach,
using absolute normality here and there. But this would be a purely
technical
trick (interesting in its own right), and we look here more for a
conceptual
analysis,
and, as long as we do not find good conceptual reasons why to be absolute
here
and not there, we will just be relative everywhere.)

(1) If $Y \xcs Z$ is a $Z-$BIG subset of $Y$ (by $Y \xcp Z),$ and $X \xcs
Y$ is a $Y-$big subset of $X$
(by $X \xcp Y),$ then $X \xcs Y \xcs Z$ is a $Z-$big subset of $X \xcs Y.$

(2) If $X \xcs Y \xcs Z$ is a $Z-$big subset of $X \xcs Y,$ and $X \xcs Y$
is a $Y-$big subset of $X$
(by $X \xcp Y)$ again, then $X \xcs Z$ is a $Z-$big subset of $X,$ so $X
\Xl  \xcp Z.$

Note that (1) is very different from Cumulativity or even Rational
Monotony, as
we do not say anything about $X$ in comparison to $Y:$ $X$ need not be any
big or
medium size subset of $Y.$

Seen as strict rules, this will not work, as it is transitivity, and thus
monotony: we have to admit exceptions, as there might just be a negative
arrow
$X \xcP Z$ in the diagram. We will discuss such situations below in (C),
where we
will modify our approach slightly, and obtain a clean analysis.

We try now to give justifications for the two (defeasible) rules. They
will be
philosophical and can certainly be contested and/or improved.

For (1):

We look at $Y.$ By $X \xcp Y,$ Y's information is accessible at $X,$ so,
as $Z-$BIG is
defined for $Y,$ $Z-$big will be defined for $Y \xcs X.$ Moreover,
there is a priori nothing
which prevents $X$ from being independent from $Y,$ i.e. $Y \xcs X$ to
behave like $Y$
with respect to $Z$ - by default: of course, there could be a negative
arrow
$X \xcP Z,$ which would prevent this.

Thus, as $Y \xcs Z$ is a $Z-$BIG subset of $Y,$ $Y \xcs X \xcs Z$ should
be a $Z-$big subset of
$Y \xcs X.$ By the same argument (independence), we should also conclude
that
$(Y-X) \xcs Z$ is a $Z-$big subset of Y-X. The definition of $Z-$big for
Y-X seems,
however, less clear.

To summarize, $Y \xcs X$ and Y-X behave by default with respect to $Z$ as
$Y$ does, i.e.
$Y \xcs X \xcs Z$ is a
$Z-$big subset of $Y \xcs X$ and $(Y-X) \xcs Z$ is a $Z-$big subset of
Y-X. The reasoning is
downward, from supersets to subsets, and symmetrical to $Y \xcs X$ and
Y-X. If the
default is violated, we need a reason for it.

This default is an assumption about the adequacy of the language. Things
do not
change wildly from one concept to another (or, better: from $Y$ to $Y \xcu
X),$ they
might change, but then we are told so - by a corresponding negative link
in the
case of diagrams.

For (2):

By $X \xcp Y,$ $X$ and $Y$ are related, and we assume that $X$ behaves as
$Y \xcs X$ does
with respect to $Z.$ This is upward reasoning, from subset to superset and
it is
NOT symmetrical: There is no reason to suppose that X-Y behaves the same
way as
$X$ or $Y \xcs X$ do with respect to $Z,$ as the only reason for $Z$ we
have, $Y,$ does not
apply.

Note that, putting relativity aside (which can also be considered as being
big/small in various, per default independent dimensions) this is close to
the
reasoning with absolutely big/small sets: $X \xcs Y-(X \xcs Y \xcs Z)$ is
small in $X \xcs Y,$ so a
fortiori small in $X,$ and $X-(X \xcs Y)$ is small in $X,$ so
$(X-(X \xcs Y)) \xcv (X \xcs Y-(X \xcs Y \xcs Z))$ is small in $X$ by the
filter property, so $X \xcs Y \xcs Z$
is big in $X,$ so a fortiori $X \xcs Z$ is big in $X.$

Thus, in summary, we conclude by default that,

(3) If $Y \xcs Z$ is a $Z-$BIG subset of $Y,$ and $X \xcs Y$ is a $Y-$big
subset of $X,$ then
$X \xcs Z$ is a $Z-$big subset of $X.$

(B) The case with longer valid paths, but without conflicts.

Treatment of longer paths: Suppose we have a valid composed path from $X$
to
$Y,$ $X \Xl  \xcp Y,$ and not any longer a direct link $X \xcp Y.$ By
induction (upward
chaining!) we argue - use directly (3) - that $X \xcs Y$ is a $Y-$big
subset of $X,$ and
conclude by (3) again that $X \xcs Z$ is a $Z-$big subset of $X.$

(C) Treatment of multiple and perhaps conflicting information.

Consider the following Diagram \ref{Diagram 5.1}:

\bfc

\bdia

\label{Diagram 5.1}

\unitlength1.0mm
\begin{picture}(130,100)

\newsavebox{\sets}
\savebox{\sets}(140,90)[bl]
{



\put(37,8){\vector(-1,1){22}}

\put(13,38){\vector(1,1){24}}
\put(67,38){\vector(-1,1){24}}

\put(40,37){\vector(0,1){23}}

\put(40,8){\vector(0,1){22}}

\put(36,34){\vector(-1,0){22}}

\put(36,4){\vector(-1,0){22}}

\put(39,3){$X$}
\put(9,33){$Y$}
\put(39,33){$Y'$}
\put(69,33){$Y''$}
\put(39,63){$Z$}
\put(9,3){$U$}

\put(38,50){\line(1,0){3.7}}

}

\put(0,0){\usebox{\sets}}
\end{picture}

\edia

\efc

We want to analyze the situation and argue that e.g. $X$ is mostly not in
$Z,$ etc.

First, all arguments about $X$ and $Z$ go via the $Y' $s. The arrows from
$X$ to
the $Y' $s, and from $Y' $ to $Y$ could also be valid paths. We look at
information
which concerns $Z$ (thus $U$ is not considered), and which is accessible
(thus $Y'' $
is not considered). (We can
be slightly more general, and consider all possible combinations of
accessible
information, not only those used in the diagram by $X.)$ Instead of
arguing
on the level of $X,$ we will argue one level above, on the Y's and their
intersections, respecting specificity and unresolved conflicts.

(Note that in more general situations, with arbitrary information
appended, the
problem is more complicated, as we have to check which information is
relevant
for some $ \xbf $ - conclusions can be arrived at by complicated means,
just as in
ordinary logic. In such cases, it might be better first to look at all
accessible information for a fixed $X,$ then at the truth values and their
relation, and calculate closure of the remaining information.)

We then have (using the obvious language: "most $ \xCB ' $s are $B' $s"
for $'' A \xcs B$ is a
big subset of $ \xCB '',$ and "MOST $ \xCB ' $s are $B' $s" for $'' A
\xcs B$ is a BIG subset of $ \xCB '' ):$

In $Y,$ $Y'',$ and $Y \xcs Y'',$ we have that MOST cases are in $Z.$
In $Y' $ and $Y \xcs Y',$ we have that MOST cases are not in $Z$ $(=$ are
in $ \xdC Z).$
In $Y' \xcs Y'' $ and $Y \xcs Y' \xcs Y'',$ we are UNDECIDED about $Z.$

Thus:

$Y \xcs Z$ will be a $Z-$BIG subset of $Y,$ $Y'' \xcs Z$ will be a $Z-$BIG
subset of $Y'',$
$Y \xcs Y'' \xcs Z$ will be a $Z-$BIG subset of $Y \xcs Y''.$

$Y' \xcs \xdC Z$ will be a $Z-$BIG subset of $Y',$ $Y \xcs Y' \xcs \xdC
Z$ will be a $Z-$BIG subset of
$Y \xcs Y'.$

$Y' \xcs Y'' \xcs Z$ will be a $Z-$MEDIUM subset of $Y' \xcs Y'',$ $Y
\xcs Y' \xcs Y'' \xcs Z$ will be
a $Z-$MEDIUM subset of $Y \xcs Y' \xcs Y''.$

This is just simple arithmetic of truth values, using specificity and
unresolved
conflicts, and the non-monotonicity is pushed into the fact that subsets
need
not preserve the properties of supersets.

In more complicated situations, we implement e.g. the general principle
(P2.2)
from Definition \ref{Definition 2.1}, to calculate the truth values. This
will use in our case specificity for conflict resolution, but it is an
abstract procedure, based on an arbitrary relation $<.$

This will result in the "correct" truth value for the intersections, i.e.
the one corresponding to the other approaches.

It remains to do two things: (C.1) We have to assure that $X$ "sees" the
correct
information, i.e. the correct intersection, and, (C.2), that $X$ "sees"
the
accepted $Y' $s, i.e. those through which valid paths go, in order to
construct not
only the result, but also the correct paths.

(Note that by split validity preclusion, if there is valid path from $
\xCB $ through $B$
to $C,$ $ \xbs:A \xFB \xcp B,$ $B \xcp C,$ and $ \xbs ':A \xFB \xcp B$
is another valid path from $ \xCB $ to $B,$ then
$ \xbs ' \xDM B \xcp C$ will also be a valid path. Proof: If not, then $
\xbs ' \xDM B \xcp C$ is precluded,
but the same preclusion will also preclude $ \xbs \xDM B \xcp C$ by split
validity
preclusion, or it is contradicted, and a similar argument applies again.
This is the argument for the simplified definition of preclusion - see
Remark \ref{Remark 2.1}, (4).)

(C.1) Finding and inheriting the correct information:

$X$ has access to $Z-$information from $Y$ and $Y',$ so we
have to consider them. Most of $X$ is in $Y,$ most of $X$ is in $Y',$
i.e.
$X \xcs Y$ is a $Y-$big subset of $X,$ $X \xcs Y' $ is a $Y' -$big subset
of $X,$ so
$X \xcs Y \xcs Y' $ is a $Y \xcs Y' -$big subset of $X,$ thus most of $X$
is in $Y \xcs Y'.$

We thus have $Y,$ $Y',$ and $Y \xcs Y' $ as possible reference classes,
and use
specificity to choose $Y \xcs Y' $ as reference class. We do not know
anything e.g. about $Y \xcs Y' \xcs Y'',$ so this is not a possible
reference class.

Thus, we use specificity twice, on the $Y' s-$level (to decide that $Y
\xcs Y' $ is
mostly not in $Z),$ and on $X' s-$level (the choice of the reference
class), but
this is good policy, as, after all, much of nonmonotonicity is about
specificity.

We should emphasize that nonmonotonicity lies in the behaviour of the
subsets
(determined by truth values and comparisons thereof) and the choice of the
reference class by specificity. But both are straightforward now and local
procedures, using information already decided before. There is no
complicated
issue here like determining extensions etc.

We now use above argument, described in the simple case, but with more
detail,
speaking in particular about the most specific reference class for
information about $Z,$ $Y \xcs Y' $ in our example - this is used
essentially in
(1.4), where the "real" information transfer happens, and here we also go
from BIG to big.

(1.1) By $X \xcp Y$ and $X \xcp Y' $ (and there are no other $Z-$relevant
information
sources), we have to consider $Y \xcs Y' $ as reference class.

(1.2) $X \xcs Y$ is a $Y-$big subset of $X$ (by $X \xcp Y)$
(it is even Y-BIG, but we are immediately more general to treat valid
paths),
$X \xcs Y' $ is a $Y' -$big subset of $X$ (by $X \xcp Y' ).$
So $X \xcs Y \xcs Y' $ is a $Y \xcs Y' -$big subset of $X.$

(1.3) $Y \xcs Z$ is a $Z-$BIG subset of $Y$ (by $Y \xcp Z),$
$Y' \xcs \xdC Z$ is a $Z-$BIG subset of $Y' $ (by $Y' \xcP Z),$
so by preclusion
$Y \xcs Y' \xcs \xdC Z$ is a $Z-$BIG subset of $Y \xcs Y'.$

(1.4) $Y \xcs Y' \xcs \xdC Z$ is a $Z-$BIG subset of $Y \xcs Y',$ and
$X \xcs Y \xcs Y' $ is a $Y \xcs Y' -$big subset of $X,$ so
$X \xcs Y \xcs Y' \xcs \xdC Z$ is a $Z-$big subset of $X \xcs Y \xcs Y'.$

This cannot be a strict rule without the reference class, as it would then
apply
to $Y \xcs Z,$ too, leading to a contradiction.

(2) If $X \xcs Y \xcs Y' \xcs \xdC Z$ is a $Z-$big subset of $X \xcs Y
\xcs Y',$ and
$X \xcs Y \xcs Y' $ is a $Y \xcs Y' -$big subset of $X,$ so
$X \xcs \xdC Z$ is a $Z-$big subset of $X.$

We make this now more formal.

We define for all nodes $X,$ $Y$ two sets: $B(X,Y),$ and $b(X,Y),$ where
$B(X,Y)$ is the set of $Y-$BIG subsets of $X,$ and
$b(X,Y)$ is the set of $Y-$big subsets of $X.$
(To distinguish undefined from medium/MEDIUM-size, we will also have to
define $M(X,Y)$ and $m(X,Y),$ but we omit this here for simplicity.)

The translations are then:

$(1.2' )$ $X \xcs Y \xbe b(X,Y)$ and $X \xcs Y' \xbe b(X,Y' )$ $ \xch $ $X
\xcs Y \xcs Y' \xbe b(X,Y \xcs Y' )$

$(1.3' )$ $Y \xcs Z \xbe B(Y,Z)$ and $Y' \xcs \xdC Z \xbe B(Y',Z)$ $ \xch
$ $Y \xcs Y' \xcs \xdC Z \xbe B(Y \xcs Y',Z)$ by preclusion

$(1.4' )$ $Y \xcs Y' \xcs \xdC Z \xbe B(Y \xcs Y',Z)$ and $X \xcs Y \xcs
Y' \xbe b(X,Y \xcs Y' )$ $ \xch $ $X \xcs Y \xcs Y' \xcs \xdC Z \xbe b(X
\xcs Y \xcs Y',Z)$
as $Y \xcs Y' $ is the most specific reference class

$(2' )$ $X \xcs Y \xcs Y' \xcs \xdC Z \xbe b(X \xcs Y \xcs Y',Z)$ and $X
\xcs Y \xcs Y' \xbe b(X,Y \xcs Y' )$ $ \xch $ $X \xcs \xdC Z \xbe b(X,Z).$

Finally:

$(3' )$ $A \xbe B(X,Y)$ $ \xcp $ $A \xbe b(X,Y)$ etc.

Note that we used, in addition to the set rules, preclusion, and the
correct
choice of the reference class.

(C.2) Finding the correct paths:

Idea:

(1) If we come to no conclusion, then no path is valid, this is trivial.

(2) If we have a conclusion:

(2.1) All contradictory paths are out: e.g. $Y \xcs Z$ will be Z-big, but
$Y \xcs Y' \xcs \xdC Z$ will be Z-big. So there is no valid path via $Y.$

(2.2) Thus, not all paths supporting the same conclusion are valid.

Consider the following Diagram \ref{Diagram 5.2}:

\bfc

\bdia

\label{Diagram 5.2}

\unitlength1.0mm
\begin{picture}(130,100)

\newsavebox{\setsx}
\savebox{\setsx}(140,90)[bl]
{



\put(37,8){\vector(-1,1){22}}
\put(41,8){\vector(1,1){22}}

\put(13,38){\vector(1,1){24}}
\put(67,38){\vector(-1,1){24}}

\put(40,37){\vector(0,1){23}}

\put(40,8){\vector(0,1){22}}

\put(36,34){\vector(-1,0){22}}
\put(66,34){\vector(-1,0){22}}


\put(39,3){$X$}
\put(9,33){$Y$}
\put(39,33){$Y'$}
\put(69,33){$Y''$}
\put(39,63){$Z$}

\put(38,50){\line(1,0){3.7}}

}

\put(0,0){\usebox{\setsx}}
\end{picture}

\edia

\efc

There might be a positive path
through $Y,$ a negative one through $Y',$ a positive one through $Y'' $
again, with
$Y'' \xcp Y' \xcp Y,$ so $Y$ will be out, and only $Y'' $ in. We can see
this, as there is a
subset, $\{Y,Y' \}$ which shows a change: $Y' \xcs Z$ is Z-BIG, $Y' \xcs
\xdC Z$ is Z-BIG,
$Y'' \xcs Z$ is Z-BIG, and $Y \xcs Y' \xcs \xdC Z$ is Z-BIG, and the
latter can only happen if
there is a preclusion between $Y' $ and $Y,$ where $Y$ looses. Thus, we
can see this
situation by looking only at the sets.

We show now equivalence with the inheritance formalism given in Section
\ref{Section 3}.

\bfa

$\hspace{0.01em}$

(+++*** Orig. No.:  Fact 5.2 )

\label{Fact 5.2}

The above definition and the one outlined in Definition \ref{Definition 2.3}
correspond.

\efa

\paragraph{
Proof:
}

$\hspace{0.01em}$

By induction on the length of the deduction that $X \xcs Z$ (or $X \xcs
\xdC Z)$ is a $Z-$big
subset of $X.$ (Outline)

It is a corollary of the proof that we have to consider only subpaths and
information of all generalized paths between $X$ and $Z.$

Make all sets (i.e. one for every node) sufficiently different, i.e.
all sets and boolean combinations of sets differ by infinitely many
elements,
e.g. $A \xcs B \xcs C$ will have infinitely many less elements than $A
\xcs B,$ etc. (Infinite
is far too many, we just choose it by laziness to have room for the
$B(X,Y)$ and
the $b(X,Y).$

Put in $X \xcs Y \xbe B(X,Y)$ for all $X \xcp Y,$ and $X \xcs \xdC Y \xbe
B(X,Y)$ for all $X \xcP Y$ as base
theory.

$Length=1:$
Then big must be BIG, and, if $X \xcs Z$ is a $Z-$BIG subset of $X,$ then
$X \xcp Z,$ likewise
for $X \xcs \xdC Z.$

We stay close now to above Diagram \ref{Diagram 5.1}, so we argue for
the negative case.

Suppose that we have deduced $X \xcs \xdC Z \xbe b(X,Z),$ we show that
there must
be a valid negative path from $X$ to $Z.$ (The other direction is easy.)

Suppose for simplicity that there is no negative link from $X$ to $Z$ -
otherwise
we are finished.

As we can distinguish intersections from elementary sets (by the starting
hypothesis about sizes), this can only be deduced using
$(2' ).$ So there must be some suitable $\{Y_{i}:i \xbe I\}$ and we must
have deduced
$X \xcs \xcS Y_{i} \xbe b(X, \xcS Y_{i}),$ the second hypothesis of $(2'
).$
If $I$ is a singleton, then we have the induction hypothesis, so there is
a valid
path from $X$ to $Y.$ So suppose $I$ is not a singleton. Then the
deduction of
$X \xcs \xcS Y_{i} \xbe b(X, \xcS Y_{i})$ can only be done by $(1.2' ),$
as this is the only rule having in the conclusion an elementary set on the
left
in $b(.,.),$ and a true intersection on the right. Going back along $(1.2'
),$
we find $X \xcs Y_{i} \xbe b(X,Y_{i}),$ and by the induction hypothesis,
there are valid paths
from $X$ to the $Y_{i}.$

The first hypothesis of $(2' ),$ $X \xcs \xcS Y_{i} \xcs \xdC Z \xbe b(X
\xcs \xcS Y_{i},Z)$ can be obtained by
$(1.3' )$ or $(1.4' ).$ If it was obtained by $(1.3' ),$ then $X$ is one
of the $Y_{i},$ but
then there is a direct link from $X$ to $Z$ (due to the "B", BIG). As a
direct link
always wins by specificity, the link must be negative, and we have a valid
negative path from $X$ to $Z.$ If it was obtained by $(1.4' ),$ then its
first
hypothesis $ \xcS Y_{i} \xcs \xdC Z \xbe B( \xcS Y_{i},Z)$ must have been
deduced, which can only be by
$(1.3' ),$ but the set of $Y_{i}$ there was chosen to take all $Y_{i}$
into account for which
there is a valid path from $X$ to $Y_{i}$ and arrows from the $Y_{i}$ to
$Z$ (the rule was
only present for the most specific reference class with respect to $X$ and
$Z!),$
and we are done by the definition of valid paths in Section \ref{Section 2}.

$ \xcz $
\\[3ex]

We summarize our ingredients.

Inheritance was done essentially by (1) and (2) of Section \ref{Section 5.1} (A)
and its elaborations (1.i), (2) and $(1.i' ),$ $(2' ).$ It consisted of a
mixture of bold and careful (in comparison to systems $P$ and $R)$
manipulation of
big subsets. We had to be bolder than the systems $P$ and $R$ are, as we
have to
transfer information also to perhaps small subsets. We had to be more
careful,
as $P$ and $R$ would have introduced far more connections than are
present. We also
saw that we are forced to loose the paradise of absolute small/big
subsets, and
have to work with relative size.

We then have a plug-in decision what to do with contradictions. This is a
plug-in, as it is one (among many possible) solutions to the much more
general
question of how to deal with contradictory information, in the presence of
a
(partial, not necessarily transitive) relation which compares strength. At
the
same place of our procedure, we can plug in other solutions, so our
approach
is truly modular in this aspect. The comparing relation is defined by the
existence of valid paths, i.e. by specificity.

This decision is inherited downward using again the specificity criterion.

Perhaps the deepest part of the analysis can be described as follows:
Relevance is coded by positive arrows, and valid positive paths, and thus
is similar to Kripke structures for modality. But, relevance (in this
reading,
which is closely related to causation) is profoundly non-monotonic, and
any
purely monotonic treatment of relevance would be insufficient. This seems
to
correspond to intuition. Relevance is then expressed formally by the
possibility
of combining different small/big sets. This is, of course, a special form
of
relevance, there might be other forms.

\section{
ACKNOWLEDGEMENTS
}

(+++*** Orig.:   ACKNOWLEDGEMENTS )
\label{Section ACKNOWLEDGEMENTS}
Some of above reflections originated directly or indirectly from a long
email discussion and criticisms with and from A.Bochman. Indirectly,
the interest of my former PhD student J.Ben-Naim in multi-valued and
information source logic has certainly influenced me to reconsider
inheritance
networks from these perspectives. Three referees have helped with their
very constructive criticism to make this text readable.

\end{document}